\documentclass[11pt,twoside]{article}
\usepackage{multicol}
\usepackage{listings}

\usepackage{amsbsy,amsfonts,amsmath,amssymb,eucal,mathrsfs}
\usepackage[all]{xy}
\usepackage{pstricks,epsfig}
\usepackage{pst-plot}
\usepackage{psfrag}
\usepackage{enumerate}
\newtheorem{definition}{Definition}[section]
\newenvironment{defi}{\begin{definition} \rm}{\end{definition}}

\newtheorem{prop}[definition]{Proposition}

\newtheorem{lemm}[definition]{Lemma}
\newtheorem{fact}[definition]{Fact}

\newtheorem{coro}[definition]{Corollary}
\newtheorem{theo}[definition]{Theorem}
\newtheorem{notation}[definition]{Notation}

\newtheorem{construction}[definition]{Construction}

\newtheorem{remark}[definition]{Remark}
\newenvironment{rema}{\begin{remark} \rm}{\end{remark}}
\newtheorem{remarks}[definition]{Remarks}

\newtheorem{example}[definition]{Example}

\newtheorem{examples}[definition]{Examples}

\newtheorem{nothing}[definition]{$\!\!$}

\newenvironment{proo}{{\flushleft \it Proof.}}{\hfill $\square$
  \vspace{2mm}}
\newenvironment{proo-prop}{{\flushleft \it Proof of Proposition
    \ref{prop-f_4}.}}{\hfill $\square$ \vspace{2mm}}

\newtheorem{conjecture}[definition]{Conjecture}

\newtheorem{definition*}{Definition}[section]
\newenvironment{defi*}{\begin{definition*} \rm}{\end{definition*}}
\newtheorem{definitions*}[definition*]{Definitions}
\newenvironment{defis*}{\begin{definitions*} \rm}{\end{definitions*}}
\newtheorem{prop*}[definition*]{Proposition}
\newtheorem{lemm*}[definition*]{Lemma}
\newtheorem{coro*}[definition*]{Corollary}
\newtheorem{theo*}[definition*]{Theorem}
\newtheorem{remark*}[definition*]{Remark}
\newenvironment{rema*}{\begin{remark*} \rm}{\end{remark*}}
\newtheorem{remarks*}[definition*]{Remarks}
\newenvironment{remas*}{\begin{remarks*} \rm}{\end{remarks*}}
\newtheorem{example*}[definition*]{Example}
\newenvironment{exam*}{\begin{example*} \rm}{\end{example*}}
\newtheorem{examples*}[definition*]{Examples}
\newenvironment{exams*}{\begin{examples*} \rm}{\end{examples*}}
\newtheorem{nothing*}[definition*]{$\!\!$}
\newenvironment{noth*}{\begin{nothing*} \rm}{\end{nothing*}}

\newtheorem{commentaire*}[definition*]{Commentaire}

\begin{document}

\def \JK {{\mathfrak{J}}}
\def \IK {{\mathfrak{I}}}
\def \Rh {{\widehat{R}}}
\def \Sh {{\widehat{S}}}
\def \supp {{\rm Supp}}
\def \codim {{\rm Codim}}
\def \b {{\beta}}
\def \T {{\Theta}}
\def \t {{\theta}}
\def \sca #1#2{\langle #1,#2 \rangle}
\def\pt{\{{\rm pt}\}}
\def\x {{\underline{x}}}
\def\y {{\underline{y}}}
\def\aut{{\rm Aut}}
\def\ra{\rightarrow}
\def\s{\sigma}\def\OO{\mathbb O}\def\PP{\mathbb P}\def\QQ{\mathbb Q}
 \def\CC{\mathbb C} \def\ZZ{\mathbb Z}\def\JO{{\mathcal J}_3(\OO)}
\newcommand{\G}{\mathbb{G}}
\def\proof{\noindent {\it Proof.}\;}
\def\qed{\hfill $\square$}
\def \uh {{\widehat{u}}}
\def \vh {{\widehat{v}}}
\def \fh {{\widehat{f}}}
\def \wh {{\widehat{w}}}
\def \Wh {{{W_{{\rm aff}}}}}
\def \Wt {{\widetilde{W}_{{\rm aff}}}}
\def \Qt {{\widetilde{Q}}}
\def \Waff {{W_{{\rm aff}}}}
\def \Waffm {{W_{{\rm aff}}^-}}
\def \Wpaff {{{(W^P)}_{{\rm aff}}}}
\def \Wtpaff {{{(\widetilde{W}^P)}_{{\rm aff}}}}
\def \Wtaffm {{\widetilde{W}_{{\rm aff}}^-}}
\def \lh {{\widehat{\lambda}}}
\def \pit {{\widetilde{\pi}}}
\def \lt {{{\lambda}}}
\def \xh {{\widehat{x}}}
\def \yh {{\widehat{y}}}
\def \a {\alpha}
\def \fll {{\longrightarrow}}
\def \b {\beta}
\def \l {\lambda}
\def \t {\theta}

\def \app {{ICI}}

%%%%%%%%% Commandes pe %%%%%%%%%%%%%%

\newcommand{\expxy}{\exp_{x \to y}}
\newcommand{\drat}{d_{\rm rat}}
\newcommand{\dmax}{d_{\rm max}}
\newcommand{\zl}{Z(x,L_x,y,L_y)}

% alg{\~A}{\AA}{!`}bres norm{\~A}{\copyright}es et anneaux usuels

\newcommand{\tiff}{if and only if }
\newcommand{\N}{\mathbb{N}}
\newcommand{\A}{{\mathbb{A}_{\rm Aff}}}
\newcommand{\Ah}{{\mathbb{A}_{\rm Aff}}}
\newcommand{\At}{{\widetilde{\mathbb{A}}_{\rm Aff}}}
\newcommand{\Ht}{{{H}^T_*(\Omega K^{\ad})}}
\renewcommand{\H}{{\rm Hi}}
\newcommand{\Ih}{{I_{\rm Aff}}}
\newcommand{\psit}{{\widetilde{\psi}}}
\newcommand{\xit}{{\widetilde{\xi}}}
\newcommand{\Jt}{{\widetilde{J}}}
\newcommand{\Zt}{{\widetilde{Z}}}
\newcommand{\Xt}{{\widetilde{X}}}
\newcommand{\at}{{\widetilde{A}}}
\newcommand{\Z}{\mathbb Z}
\newcommand{\R}{\mathbb{R}}
\newcommand{\Q}{\mathbb{Q}}
\newcommand{\C}{\mathbb{C}}
\renewcommand{\O}{\mathbb{O}}
\newcommand{\F}{\mathbb{F}}
\newcommand{\p}{\mathbb{P}}
\newcommand{\co}{{\cal O}}
\newcommand{\pos}{{\bf P}}

\renewcommand{\a}{{\alpha}}
\newcommand{\az}{\a_\Z}
\newcommand{\ak}{\a_k}

\newcommand{\rc}{\R_\C}
\newcommand{\cc}{\C_\C}
\newcommand{\hc}{\H_\C}
\newcommand{\oc}{\O_\C}

\newcommand{\rk}{\R_k}
\newcommand{\ck}{\C_k}
\newcommand{\hk}{\H_k}
\newcommand{\ok}{\O_k}

\newcommand{\rz}{\R_Z}
\newcommand{\cz}{\C_Z}
\newcommand{\hz}{\H_Z}
\newcommand{\oz}{\O_Z}

\newcommand{\RR}{\R_R}
\newcommand{\CR}{\C_R}
\newcommand{\HR}{\H_R}
\newcommand{\OR}{\O_R}

\newcommand{\re}{\mathtt{Re}}

\newcommand{\matttr}[9]{
\left (
\begin{array}{ccc}
{} \hspace{-.2cm} #1 & {} \hspace{-.2cm} #2 & {} \hspace{-.2cm} #3 \\
{} \hspace{-.2cm} #4 & {} \hspace{-.2cm} #5 & {} \hspace{-.2cm} #6 \\
{} \hspace{-.2cm} #7 & {} \hspace{-.2cm} #8 & {} \hspace{-.2cm} #9
\end{array}
\hspace{-.15cm}
\right )   }

% alg{\~A}{\AA}{!`}bre

\newcommand{\dual}{{\bf v}}
\newcommand{\com}{\mathtt{Com}}
\newcommand{\rg}{\mathtt{rg}}
\newcommand{\pu}{{\mathbb{P}^1}}
\newcommand{\scal}[1]{\langle #1 \rangle}
\newcommand{\MK}[2]{{\overline{{\rm M}}_{#1}(#2)}}
\newcommand{\MKE}[3]{{{{\rm M}}_{#1}^{#2}(#3)}}
\newcommand{\SM}[2]{{\overline{{\rm\bf M}}_{#1}(#2)}}
\newcommand{\SME}[3]{{\overline{{\rm\bf M}}_{#1}^{#2}(#3)}}
\newcommand{\mor}[2]{{{\rm Mor}_{#1}(\pu,#2)}}

\newcommand{\fg}{\mathfrak g}
\newcommand{\fgad}{{\mathfrak g}^{\rm ad}}
\renewcommand{\fh}{\mathfrak h}
\newcommand{\fu}{\mathfrak u}
\newcommand{\fz}{\mathfrak z}
\newcommand{\fn}{\mathfrak n}
\newcommand{\fe}{\mathfrak e}
\newcommand{\fp}{\mathfrak p}
\newcommand{\ft}{\mathfrak t}
\newcommand{\fl}{\mathfrak l}
\newcommand{\fq}{\mathfrak q}
\newcommand{\fsl}{\mathfrak {sl}}
\newcommand{\fgl}{\mathfrak {gl}}
\newcommand{\fso}{\mathfrak {so}}
\newcommand{\fsp}{\mathfrak {sp}}
\newcommand{\ff}{\mathfrak {f}}

\newcommand{\ad}{{\rm ad}}
\newcommand{\jad}{{j^\ad}}
\newcommand{\id}{{\rm id}}

%%%% Poids et racines   %%%%

\newcommand{\dynkinadeux}[2]
{
$#1$
\setlength{\unitlength}{1.2pt}
\hspace{-3mm}
\begin{picture}(12,3)
\put(0,3){\line(1,0){10}}
\end{picture}
\hspace{-2.4mm}
$#2$
}

\newcommand{\mdynkinadeux}[2]
{
\mbox{\dynkinadeux{#1}{#2}}
}

\newcommand{\dynkingdeux}[2]
{
$#1$
\setlength{\unitlength}{1.2pt}
\hspace{-3mm}
\begin{picture}(12,3)
\put(1,.8){$<$}
\multiput(0,1.5)(0,1.5){3}{\line(1,0){10}}
\end{picture}
\hspace{-2.4mm}
$#2$
}

\newcommand{\poidsesix}[6]
{
\hspace{-.12cm}
\left (
\begin{array}{ccccc}
{} \hspace{-.2cm} #1 & {} \hspace{-.3cm} #2 & {} \hspace{-.3cm} #3 &
{} \hspace{-.3cm} #4 & {} \hspace{-.3cm} #5 \vspace{-.13cm}\\
\hspace{-.2cm} & \hspace{-.3cm} & {} \hspace{-.3cm} #6 &
{} \hspace{-.3cm} & {} \hspace{-.3cm}
\end{array}
\hspace{-.2cm}
\right )      }

\newcommand{\copoidsesix}[6]{
\hspace{-.12cm}
\left |
\begin{array}{ccccc}
{} \hspace{-.2cm} #1 & {} \hspace{-.3cm} #2 & {} \hspace{-.3cm} #3 &
{} \hspace{-.3cm} #4 & {} \hspace{-.3cm} #5 \vspace{-.13cm}\\
\hspace{-.2cm} & \hspace{-.3cm} & {} \hspace{-.3cm} #6 &
{} \hspace{-.3cm} & {} \hspace{-.3cm}
\end{array}
\hspace{-.2cm}
\right |      }

\newcommand{\poidsesept}[7]{
\hspace{-.12cm}
\left (
\begin{array}{cccccc}
{} \hspace{-.2cm} #1 & {} \hspace{-.3cm} #2 & {} \hspace{-.3cm} #3 &
{} \hspace{-.3cm} #4 & {} \hspace{-.3cm} #5 & {} \hspace{-.3cm} #6
\vspace{-.13cm}\\
\hspace{-.2cm} & \hspace{-.3cm} & {} \hspace{-.3cm} #7 &
{} \hspace{-.3cm} & {} \hspace{-.3cm}
\end{array}
\hspace{-.2cm}
\right )      }

\newcommand{\copoidsesept}[7]{
\hspace{-.12cm}
\left |
\begin{array}{cccccc}
{} \hspace{-.2cm} #1 & {} \hspace{-.3cm} #2 & {} \hspace{-.3cm} #3 &
{} \hspace{-.3cm} #4 & {} \hspace{-.3cm} #5 & {} \hspace{-.3cm} #6
\vspace{-.13cm}\\
\hspace{-.2cm} & \hspace{-.3cm} & {} \hspace{-.3cm} #7 &
{} \hspace{-.3cm} & {} \hspace{-.3cm}
\end{array}
\hspace{-.2cm}
\right |      }

\newcommand{\poidsehuit}[8]{
\hspace{-.12cm}
\left (
\begin{array}{cccccc}
{} \hspace{-.2cm} #1 & {} \hspace{-.3cm} #2 & {} \hspace{-.3cm} #3 &
{} \hspace{-.3cm} #4 & {} \hspace{-.3cm} #5 & {} \hspace{-.3cm} #6 &
{} \hspace{-.3cm} #7   \vspace{-.13cm}\\
\hspace{-.2cm} & \hspace{-.3cm} & {} \hspace{-.3cm} #8 &
{} \hspace{-.3cm} & {} \hspace{-.3cm}
\end{array}
\hspace{-.2cm}
\right )      }

\newcommand{\copoidsehuit}[8]{
\hspace{-.12cm}
\left |
\begin{array}{cccccc}
{} \hspace{-.2cm} #1 & {} \hspace{-.3cm} #2 & {} \hspace{-.3cm} #3 &
{} \hspace{-.3cm} #4 & {} \hspace{-.3cm} #5 & {} \hspace{-.3cm} #6 &
{} \hspace{-.3cm} #7  \vspace{-.13cm}\\
\hspace{-.2cm} & \hspace{-.3cm} & {} \hspace{-.3cm} #8 &
{} \hspace{-.3cm} & {} \hspace{-.3cm}
\end{array}
\hspace{-.2cm}
\right |      }

\newcommand{\im}{\mathtt{Im}}

%%%% Caligraphic letters %%%%%%%%%%%%%%

\def\cA{{\cal A}} \def\cC{{\cal C}} \def\cD{{\cal D}} \def\cE{{\cal E}}
\def\cF{{\cal F}} \def\cG{{\cal G}} \def\cH{{\cal H}} \def\cI{{\cal I}}
\def\cK{{\cal K}} \def\cL{{\cal L}} \def\cM{{\cal M}} \def\cN{{\cal N}}
\def\cO{{\cal O}}
\def\cP{{\cal P}} \def\cQ{{\cal Q}} \def\cT{{\cal T}} \def\cU{{\cal U}}
\def\cV{{\cal V}} \def\cX{{\cal X}} \def\cY{{\cal Y}} \def\cZ{{\cal Z}}

\def \g {{\gamma}}

\def \cw {{\underline{w}}}
\def \wt {{\widetilde{w}}}
\def \ga {{\gamma}}
\def \at {{\widetilde{\alpha}}}
\def \Hom {{\rm Hom}}
\def \pr {{\rm pr}}
\def \im {{\rm Im}}
\def \pic {{\rm Pic}}
\def \Stab {{\rm Stab}}
\def \we {{\rm Weil}}
\def \pu {{\mathbb P}^1}
\def \opu {{\oo_\pu}}
\def \Mor #1#2{{\bf{Hom}}_{#1}(\pu,#2)}
\def \MorC #1#2{{\bf{Hom}}_{#1}(C,#2)}
\def \ne {N\!E}
\def \nef {{\rm Nef}}
\def \neff {{\rm Neff}}
\def \vp {{\phi}}

\def \pp {\odot}
\def \tr {{}^t}
\def \ct {{}_Tc}
\def \lcom {$\Lambda$-(co)minuscule }
\def \lmin {$\Lambda$-minuscule }
\renewcommand{\ok}[2]{{#1} \cdot {#2} = {#1} \pp {#2}}
\newcommand{\notok}[2]{{#1} \cdot {#2} \not = {#1} \pp {#2}}
\newcommand{\oknu}[3]{c_{{#1},{#2}}^{{#3}} = t_{{#1},{#2}}^{{#3}} \cdot
m_{{#1},{#2}}^{{#3}}}
\newcommand{\notoknu}[3]{c_{{#1},{#2}}^{{#3}} \not = t_{{#1},{#2}}^{{#3}} \cdot
m_{{#1},{#2}}^{{#3}}}

 \title{Elliptic curves on some homogeneous spaces}
 \author{B. Pasquier \& N. Perrin}

\maketitle

\begin{abstract}
Let $X$ be a minuscule homogeneous space, an odd quadric, or an adjoint homogenous space of type different from $A$ and $G_2$. Le $C$ be an elliptic curve. In this paper, we prove that for $d$ large enough, the scheme of degree $d$ morphisms from $C$ to $X$ is irreducible, giving an explicit lower bound for $d$ which is optimal in many cases.
\end{abstract}

 {\def\thefootnote{\relax}
 \footnote{ \hspace{-6.8mm}
 }
 }

%%%%%%%%%%%%%%%%%%%%%%%%% INTRO %%%%%%%%%%%%%%%%%%%%%%%%%%%%%%%%%%

\section*{Introduction}

In this paper we study the scheme $\MorC{d}{X}$ of degree $d$ morphisms from an elliptic curve $C$ to a rational homogeneous space $X$. Specifically we will assume that $X$ is either a minuscule homogeneous space or an adjoint homogeneous space. Minuscule homogeneous spaces are natural generalisation of Grassmann varieties (see Table 1). Adjoint homogeneous space are also called quasi-minuscule and are obtained as the unique closed orbits in $\p\mathfrak{g}$ under the adjoint action. These varieties are also called minimal nilpotent orbit (see Table 2).

We prove that the scheme $\MorC{d}{X}$ is irreducible as soon as $d$ is large enough and we give an explicit bound $d(X)$ for $d$ (see Tables 1 and 2).

\begin{theo}
\label{main}
Let $C$ be a smooth elliptic curve and let $X$ be a minuscule or an adjoint homogeneous space. If $X$ is of adjoint type, assume furthermore that the group $G$ is not of type $A$ or $G_2$. For $d\geq d(X)$, the scheme of morphisms $\MorC{d}{X}$ is irreducible of dimension $c_1(X)d$.
\end{theo}

In many cases, we also prove that the bound $d(X)$ is optimal in the following sense: for $d<d(X)$, the scheme $\MorC{d}{X}$ has dimension strictly bigger than the expected dimension $c_1(X)d$. Note that it may happen that for $d<d(X)$ the scheme $\MorC{d}{X}$ is irreducible. This is the case for all $d$ when $X$ is isomorphic to the maximal isotropic Grassmann variety $\G_Q(n,2n)$ (see Table 1 and Corollary \ref{coro-gq}).

\begin{rema}
(\i) Note that B. Kim and R. Pandharipande proved the connectedness of the moduli space of stable maps $\overline{{\bf M}}_{d,g}(X)$ to any homogeneous space and for any genus. However, this space is almost never irreducible since there will be irreducible components of unexpected dimension (see Proposition \ref{non-irred}).

(\i\i) Note also that for $d=2$, the scheme of morphisms $\MorC{d}{X}$ to any homogeneous space with Picard rank $1$ is irreducible of dimension $c_1(X)+\dim X+1$. Indeed, any degree $2$ morphism factors through a line.
\end{rema}

Minuscule and adjoint homogeneous spaces are of the form $G/P$ for $G$ a semisimple algebraic group. Except for the adjoint homogeneous space in type $A$, the parabolic subgroup $P$ is a maximal parabolic subgroup. This is the reason why we will assume that if $X$ is of adjoint type the group $G$ is not of type $A$. Recall that the Dynkin diagram of the group $G$ has vertices indexed by simple roots.
In the following tables, we give the list of all minuscule homogeneous spaces. Note that we also include odd dimensional quadrics which are not minuscule but cominuscule. The above statement is still true for odd dimensional quadrics.

\begin{equation*}
\label{table}
\begin{array}{ccccccc}
\hline
Type & Variety & Diagram  & Dimension &
\hspace*{5mm}Index \hspace*{5mm} & d(X) \\
\hline
 A_{n-1} &\hspace{5mm} \G(k,n)\hspace{5mm} &\setlength{\unitlength}{2.5mm}
\begin{picture}(15,3)(-2,0)
\put(0,0){$\circ$}
\multiput(2,0)(2,0){5}{$\circ$}
\multiput(.73,.4)(2,0){5}{\line(1,0){1.34}}
\put(4,0){$\bullet$}
\end{picture} & k(n-k) & n & n \\
B_n & \hspace{5mm}\Q_{2n-1}\hspace{5mm} &
\setlength{\unitlength}{2.5mm}
\begin{picture}(15,3)(-2,0)
\put(0,0){$\circ$}
\multiput(2,0)(2,0){5}{$\circ$}
\multiput(.73,.4)(2,0){4}{\line(1,0){1.34}}
\multiput(8.73,.2)(0,.4){2}{\line(1,0){1.34}}
\put(0,0){$\bullet$}
\end{picture} & 2n-1 & 2n-1 & 3\\
D_n & \hspace{5mm}\Q_{2n-2}\hspace{5mm} &
\setlength{\unitlength}{2.5mm}
\begin{picture}(15,3)(-2,0)
\put(2,0){$\circ$}
\multiput(4,0)(2,0){4}{$\circ$}
\multiput(2.73,.4)(2,0){4}{\line(1,0){1.34}}
\put(10,0){$\bullet$}
\put(0,-1.1){$\circ$}
\put(0,1.2){$\circ$}
\put(.6,1.5){\line(5,-3){1.5}}
\put(.6,-.64){\line(5,3){1.5}}
\end{picture}
\vspace{.2cm}
 & 2n-2 & 2n-2  & 3\\
D_n & \hspace{5mm}\G_Q(n,2n)\hspace{5mm} &
\setlength{\unitlength}{2.5mm}
\begin{picture}(15,3)(-2,0)
\put(2,0){$\circ$}
\multiput(4,0)(2,0){4}{$\circ$}
\multiput(2.73,.4)(2,0){4}{\line(1,0){1.34}}
\put(0,1.2){$\bullet$}
\put(0,-1.1){$\circ$}
\put(.6,1.5){\line(5,-3){1.5}}
\put(.6,-.64){\line(5,3){1.5}}
\end{picture}
\vspace{.2cm}
 & \frac{n(n-1)}{2} & 2n-2 & n-1 \\
E_6 &\hspace{5mm} \OO\PP^2 \hspace{5mm}&
\setlength{\unitlength}{2.5mm}
\begin{picture}(15,3)(-3,-.5)
\put(0,0){$\circ$}
\multiput(2,0)(2,0){4}{$\circ$}
\multiput(.73,.4)(2,0){4}{\line(1,0){1.34}}
\put(0,0){$\bullet$}
\put(4,-2){$\circ$}
\put(4.42,-1.28){\line(0,1){1.36}}
\end{picture}
& 16 & 12 & 3\\
E_7 & \hspace{5mm}E_7/P_7\hspace{5mm} &
\setlength{\unitlength}{2.5mm}
\begin{picture}(15,3)(-2,0)
\put(0,0){$\circ$}
\multiput(2,0)(2,0){5}{$\circ$}
\multiput(.73,.4)(2,0){5}{\line(1,0){1.34}}
\put(10,0){$\bullet$}
\put(4,-2){$\circ$}
\put(4.42,-1.28){\line(0,1){1.36}}
\end{picture}
& 27 & 18 & 8\\
\\
\hline
\\
\multicolumn{6}{c}{\textrm{Table 1. Minuscule homogeneous spaces and odd quadrics.}}\\
\end{array}
\end{equation*}

In the following table, we give the list of all adjoint homogeneous spaces for a group of type different from $A$. In this table we also include the adjoint variety of type $G_2$ even if we have no result for this variety.

\begin{equation*}
\begin{array}{ccccccc}
\hline
Type & Variety & Diagram  & Dimension &
\hspace*{5mm}Index \hspace*{5mm} & d(X) \\
\hline
B_n & \hspace{5mm}\G_Q(2,2n+1)\hspace{5mm} &
\setlength{\unitlength}{2.5mm}
\begin{picture}(15,3)(-2,0)
\put(0,0){$\circ$}
\multiput(2,0)(2,0){5}{$\circ$}
\multiput(.73,.4)(2,0){4}{\line(1,0){1.34}}
\multiput(8.73,.2)(0,.4){2}{\line(1,0){1.34}}
\put(2,0){$\bullet$}
\end{picture} & 4n-5 & 2n-2 & 2n\\
C_n & \hspace{5mm}\p^{2n-1}\hspace{5mm} &
\setlength{\unitlength}{2.5mm}
\begin{picture}(15,3)(-2,0)
\put(0,0){$\circ$}
\multiput(2,0)(2,0){5}{$\circ$}
\multiput(.73,.4)(2,0){4}{\line(1,0){1.34}}
\multiput(8.73,.2)(0,.4){2}{\line(1,0){1.34}}
\put(0,0){$\bullet$}
\end{picture} & 2n-1 & 2n & 2 \\
D_n & \hspace{5mm}\G_Q(2,2n)\hspace{5mm} &
\setlength{\unitlength}{2.5mm}
\begin{picture}(15,3)(-2,0)
\put(2,0){$\circ$}
\multiput(4,0)(2,0){4}{$\circ$}
\multiput(2.73,.4)(2,0){4}{\line(1,0){1.34}}
\put(8,0){$\bullet$}
\put(0,-1.1){$\circ$}
\put(0,1.2){$\circ$}
\put(.6,1.5){\line(5,-3){1.5}}
\put(.6,-.64){\line(5,3){1.5}}
\end{picture}
\vspace{.2cm}
 & 4n-7 & 2n-3  & 2n-1\\
E_6 &\hspace{5mm} E_6/P_2 \hspace{5mm}&
\setlength{\unitlength}{2.5mm}
\begin{picture}(15,3)(-3,-.5)
\put(0,0){$\circ$}
\multiput(2,0)(2,0){4}{$\circ$}
\multiput(.73,.4)(2,0){4}{\line(1,0){1.34}}
\put(4,-2){$\bullet$}
\put(4,-2){$\circ$}
\put(4.42,-1.28){\line(0,1){1.36}}
\end{picture}
& 21 & 11 & 9\\
E_7 & \hspace{5mm}E_7/P_7\hspace{5mm} &
\setlength{\unitlength}{2.5mm}
\begin{picture}(15,3)(-2,0)
\put(0,0){$\circ$}
\multiput(2,0)(2,0){5}{$\circ$}
\multiput(.73,.4)(2,0){5}{\line(1,0){1.34}}
\put(0,0){$\bullet$}
\put(4,-2){$\circ$}
\put(4.42,-1.28){\line(0,1){1.36}}
\end{picture}
& 33 & 17 & 11\\
E_8 & \hspace{5mm}E_8/P_8\hspace{5mm} &
\setlength{\unitlength}{2.5mm}
\begin{picture}(15,3)(-2,0)
\put(0,0){$\circ$}
\multiput(2,0)(2,0){5}{$\circ$}
\multiput(.73,.4)(2,0){6}{\line(1,0){1.34}}
\put(12,0){$\bullet$}
\put(4,-2){$\circ$}
\put(4.42,-1.28){\line(0,1){1.36}}
\end{picture}
& 57 & 29 & 15\\
F_4 & \hspace{5mm}F_4/P_1\hspace{5mm} &
\setlength{\unitlength}{2.5mm}
\begin{picture}(15,3)(-2,0)
\put(0,0){$\circ$}
\multiput(2,0)(2,0){1}{$\circ$}
\multiput(4,0)(2,0){1}{$\circ$}
\multiput(6,0)(2,0){1}{$\circ$}
\multiput(.73,.4)(2,0){1}{\line(1,0){1.34}}
\multiput(4.73,.4)(2,0){1}{\line(1,0){1.34}}
\multiput(2.73,.2)(0,.4){2}{\line(1,0){1.34}}
\put(0,0){$\bullet$}
\end{picture} & 15 & 8 & 8 \\
G_2 & \hspace{5mm}G_2/P_1\hspace{5mm} &
\setlength{\unitlength}{2.5mm}
\begin{picture}(15,3)(-2,0)
\put(0,0){$\circ$}
\multiput(2,0)(2,0){1}{$\circ$}
\multiput(.73,.4)(2,0){1}{\line(1,0){1.34}}
\multiput(0.73,.4)(2,0){1}{\line(1,0){1.34}}
\multiput(0.73,.2)(0,.4){2}{\line(1,0){1.34}}
\put(0,0){$\bullet$}
\end{picture} & 5 & 3 &  \\
\\
\hline
\\
\multicolumn{6}{c}{\textrm{Table 2. Adjoint homogeneous spaces of Picard number one.}}\\
\end{array}
\end{equation*}

In these tables, we followed the notation of \cite[Tables]{bou} and we depicted the set $\Sigma(P)$ of simple roots not in $P$ with plain vertices. The minuscule and adjoint varieties $X=G/P$ are described in the second column. By convention we denote by $\G(k,n)$ (resp. $\G_Q(k,n)$) the Grassmann variety of $k$-dimensional subspaces in $\C^n$ (resp. isotropic $k$-dimensional subspaces in $\C^{n}$ for a non-degenerate quadratic form $Q$ in $\C^{n}$). For $\G_Q(n,2n)$ we only consider one of the two connected components of the above Grassmann variety. We denoted by $\QQ_m$ any smooth $m$-dimensional quadric. The varieties $\OO \PP^2=E_6/P_1$ and $E_7/P_7$ are the Cayley plane and the Freudenthal variety. Recall that the index of $X=G/P$ is the integer $c_1(X)$ such that the anticanonical divisor $-K_X$ equals $c_1(X)H$ where $H$ is an ample generator of the Picard group. The last column gives the bound $d(X)$ of our main theorem.

\begin{rema}
Note that the above statement was already known for Grassmann varieties by results of A. Brugui{\`e}res \cite{bruguieres} and for quadrics by results of E. Ballico \cite{ballico}. Note also that for orthogonal Grassmann varieties, the above result was obtained by another technique by the second author in the unpublished paper \cite{arxiv}. This case was the very first motivation for our study. It answered a question of D. Markushevitch and has been used in \cite{marku}.
\end{rema}

Let us briefly describe the content of the paper. The strategy of the proof relies on a description given in \cite{fourier} of a \emph{big} open cell $U$ in any homogeneous space $X$ as a tower of affine bundles $\phi:U\to Y$ where $Y$ is again homogeneous under a smaller group. In section \ref{reduction} we explain how to restrict the study of the irreducibility of $\MorC{d}{X}$ to the one of $\MorC{d}{U}$ and we explain, under some conditions on the fibration $\phi:U\to Y$,
how to deduce the irreducibility of $\MorC{d}{U}$ from the irreducibility of $\MorC{d}{Y}$. In section \ref{sect-fib}, we recall the construction of $U$ and of the the fibration $\phi:U\to Y$ and describe the tower of affine fibrations.
In section \ref{sect-annu}, we prove some results on the cohomology of the restriction to elliptic curves of equivariant locally free sheaves on projective spaces
quadrics. In section \ref{sect-irr}, we gather all these results to prove Theorem \ref{main}. In section \ref{sect-stable}, we include some remarks on the moduli space of stable maps to minuscule or adjoint homogeneous spaces with source an elliptic curve. Finally, in an appendix we deal with of the quadric of dimension~3, which we cannot consider with the same method as all other homogeneous spaces.

\vskip 0.3 cm

{\bf Acknowledgments.} The second author would like to thank
Christian Peskine for enlightening discussions on the paper \cite{GLP} and
Piotr Achinger for useful comments on Spinor bundles.

\tableofcontents

%%%%%%%%%%%%%%%%%%%%%%%%%%% 1 %%%%%%%%%%%%%%%%%%%%%%%%

\section{Reduction of the theorem}
\label{reduction}

\subsection{Restriction to an open subset}

Let $C$ be a smooth curve, let $X$ be a homogeneous space with Picard number 1 and let $\a\in A_1(X)$. The strategy of the proof will be similar to the one in \cite{fourier}. We prove that there exists an open subset $U$ of $X$ whose complement is of codimension at least 2 and such that $U$ can be realised as a tower of affine bundles $\phi:U\to Y$ over a homogeneous space $Y$ of smaller dimension. The next result from \cite[Proposition 2]{fourier} proves that the irreducibility of $\MorC{\a}{X}$ is equivalent to the irreducibility of $\MorC{i^*\a}{U}$ (where $i:U\to X$ is the inclusion morphism). This result is an easy application of Kleiman-Bertini Theorem \cite{kleiman}.

\begin{prop}
\label{ouvert}
Let $C$ be a smooth curve and let $X$ be a homogeneous space under a group $G$. Assume that $U$ be an open subset of $X$ such that $\codim_X(X\setminus U)\geq2$ and let $i:U\to X$ be the inclusion morphism. For $\a\in A_1(X)$, if $\MorC{i^*\a}{U}$ is irreducible, then so is $\MorC{\a}{X}$.
\end{prop}

We then need to study the open subset $U$. It will be the open orbit of the action of some parabolic subgroup of $G$ on $X$. The second author proved in \cite{fourier} that $U$ can be realised as a tower $\phi:U\to Y$ of affine bundles over a homogeneous space $Y$ of smaller dimension. However, if this was enough to study rational curves, we need to be more precise if we want to study higher genus curves. If $\phi:U\to Y$ is an affine bundle associated to a vector bundles $E$, we explain
in the next subsection how a good understanding of the restriction of $E$ to higher genus curves in $Y$ enables to proceed by induction on the dimension.

\subsection{Passing through affine bundles}

Let us fix some notation. In this section $X$ is a homogeneous space of Picard number 1, we denote by $U$ an open subset whose complement $Z$ is of codimension 2 and we assume that $U$ can be realised as a sequence of affine bundles $\phi:U\to Y$ over a homogeneous space $Y$ of smaller dimension.

Note that because of the inequality $\codim_X Z\geq2$, we have $\pic(U)=\pic(X)$. Furthermore since $\phi:U\to Y$ is a sequence of affine bundles, then we also have $\pic(U)=\pic(Y)$. For $f:C\to U$ a morphism, we may define a linear form $[f]$ on $\pic(X)$ by ${\cal L}\mapsto\deg(f^*{\cal L})$. We shall say that for $\a\in\pic(X)^\vee$ a curve $f:C\to U$ is of class $\a$ if $[f]=\a$. We shall write $\scal{\a,{\cal L}}$ for the evaluation of $\a$ at ${\cal L}\in\pic(X)$.
For $\a$ such a class in $\pic(X)^\vee$, we denote by $\MorC{\a}{X}$
resp. $\MorC{\a}{U}$, $\MorC{\a}{Y}$
the scheme of morphisms from $f:C\to X$
resp. $f:C\to U$, $f:C\to Y$ such that $[f]=\a$. Note that for $X$ (or for $Y$) we have identifications $A_1(X)\simeq\pic(X)^\vee$ and $A_1(Y)\simeq\pic(Y)^\vee$
and that these schemes of morphisms are the classical schemes of morphisms as defined in \cite{grothendieck}.

We now prove a general result on the scheme of morphisms to an affine bundle. For this we introduce the following notation. Let $\phi:U\to Y$ be an affine bundle whose direction vector bundle is associated to a locally free sheaf $E$ on $Y$ and let $\a\in \pic(U)^\vee$. We denote by ${\bf H}_b$ the locally closed subset of ${\bf Hom}_{\vp_*\a}(C,Y)$ defined by
$${\bf H}_b=\{f\in{\bf Hom}_{\vp_*\a}(C,Y)\ /\ \dim H^1(C,f^*E)=b\}.$$

\begin{prop}
\label{prop-vb}
Let $\vp:U\to Y$ be an affine bundle over a variety $Y$ whose direction vector bundle is associated to a locally free sheaf $E$. Let $\a\in{\rm Pic}(U)^\vee$ and let $C$ be a smooth curve such that the scheme ${\bf Hom}_{\vp_*\a}(C,Y)$ is irreductible of dimension $\scal{\vp_*\a,c_1(T_Y)}$.

(\i) Assume that there exists $f\in{\bf Hom}_{\vp_*\a}(C,Y)$ satisfying the vanishing $H^1(C,f^*E)=0$. If for all integer $b>0$ we have the inequality $\codim\,{\bf H}_b>b$,
then ${\bf Hom}_{\a}(C,U)$ is irreducible of dimension
$\scal{\a,c_1(T_U)}$.

(\i\i) Assume that $\vp$ is a vector bundle and let $a$ be the smallest integer such that there exists an element $f\in{\bf Hom}_{\vp_*\a}(C,Y)$ with $\dim H^1(C,f^*E)=a$. If there exists an integer $b>a$ satisfying the inequality $\codim\,{\bf H}_b\leq b-a$,
then the scheme ${\bf Hom}_{\a}(C,U)$ is reductible.
\end{prop}

\begin{proo}
(\i) Let us consider the natural morphism ${\bf Hom}_{\a}(C,U)\to {\bf Hom}_{\vp_*\a}(C,Y)$.
As explained in \cite{fourier} for rational curves, the fiber of this map over an element $f\in {\bf Hom}_{\vp_*\a}(C,Y)$ is given by the sections $s\in H^0(C,f^*E)$.  In particular, such a section exists if and only if the affine bundle $C\times_YU\to C$ obtained by pull-back to $C$ is a vector bundle. For a general element $f\in{\bf Hom}_{\phi_*\a}(C,Y)$ this is the case thanks to the vanishing $H^1(C,f^*E)=0$.

For $b>0$ the fibers over the locally closed subset ${\bf H}_b$ are either empty (if the pull-back is not a vector bundle) or of
dimension $\dim H^0(C,f^*E)=\chi(C,f^*E)+\dim H^1(C,f^*E)=\scal{\phi_*\a,c_1(E)}+b$. The dimension of the inverse image of ${\bf H}_b$ in ${\bf Hom}_{\a}(C,U)$ is therefore of dimension strictly less than $\scal{\phi_*\a,c_1(T_Y)+c_1(E)}=\scal{\a,c_1(T_U)}$ which is the expected dimension of ${\bf Hom}_{\a}(C,U)$. As any irreducible component is at least of this expected dimension (cf. \cite{MO}) these locally closed subsets do not form irreducible components. Therefore the inverse image of the open subset ${\bf H}_0$
is dense in  ${\bf Hom}_{\a}(C,U)$. As explained in \cite{fourier} for rational curves, this inverse image is an affine bundle over the base and is thus irreducible of dimension  $\scal{\phi_*\a,c_1(T_Y)+c_1(E)}=\scal{\a,c_1(T_U)}$.

(\i\i) In this case, the inverse image of the locally closed subset ${\bf H}_a$
is of dimension $\scal{\a,c_1(T_U)}+a$ while the inverse image of the locally closed subset ${\bf H}_b$
is of dimension at least $\scal{\a,c_1(T_U)}+a$. This second locally closed subset is therefore not contained in the closure of the first one thus they have to be contained in two different irreducible components.
\end{proo}

%%%%%%%%%%%%%%%%%%%%%%%%%%%%%%%%%% 2 %%%%%%%%%%%%%%%%%%%%%%%%

\subsection{Decomposition of $X$}
\label{sect-fib}

The aim of this subsection is to describe a \emph{big} open cell $U$ in the homogeneous space $X$ and to realise it as a fibration $\phi:U\to Y$ over a space $Y$ homogeneous under a group of smaller rank. This fibration can be decomposed in general as a sequence of affine bundles (see \cite[Proposition 5]{fourier}). Let us recall this result.

Let $T$ be a maximal torus in $G$ and let $B$ be a Borel subgroup of $G$ containing $T$. Denote by $w_0$ the longest element of the Weyl group $W$ of $(G,T)$. For any subgroup $H$ of $G$, we set $H^{w_0}:=w_0Hw_0^{-1}$ and for any parabolic subgroup $P$ of $G$ containing $B$, we denote by $L_P$ the Levi subgroup of $P$ containing $T$ and by $U_P$ the unipotent radical of $P$. We also denote by $U_P^r\subset\cdots\subset U_P^1$ the central filtration of $U_P$. For any parabolic subgroup $P$ containing $B$, we denote by $\Sigma(P)$ the set of simple roots which are not roots of $P$. Finally, for $\a$ a simple root, we denote by $\iota(\a)$ its image via the Weyl involution \emph{i.e.} $\iota(\a)=-w_0(\a)$.

\begin{lemm}[Proposition 5 and 6, \cite{fourier}]\label{lemm:affbundle}
Let $P$ and $Q$ be two parabolic subgroup of $G$ containing $B$ and assume that $\iota(\Sigma(P))\cap \Sigma(Q)=\emptyset$. Set $U:=Q^{w_0}P/P$. Then $U$ is an open subset of $G/P$ with complement in codimension at least 2 and is an $L_{Q^{w_0}}$ equivariant fibration $\vp:U\to Y$ over the flag variety $Y=L_{Q^{w_0}}/(L_{Q^{w_0}}\cap P)$. Moreover this fibration decomposes in a sequence of affine bundles whose associated vector bundles are defined over $Y$ and are the $L_{Q^{w_0}}$-equivariant vector bundles associated to the $L_{Q^{w_0}}\cap P$ representations $U^i_{Q^{w_0}}/U^{i+1}_{Q^{w_0}}\cdot (U^i_{Q^{w_0}}\cap P)$.
\end{lemm}

\begin{rema}
Note that it is easy to determine the weights of $T$ on the representation $U^i_{Q^{w_0}}/U^{i+1}_{Q^{w_0}}\cdot (U^i_{Q^{w_0}}\cap P)$ defining the vector bundle. For example for $P$ and $Q$ maximal parabolic subgroups with $\Sigma(P)=\{\a_P\}$ and $\Sigma(Q)=\{\a_Q\}$, then the weights of the above representation are the negative roots $\a$ such that the coefficient of $\a$ on $\a_P$ is negative and the coefficient of $\a$ on $\a_Q$ is exactly $-i$.
\end{rema}

Keeping the previous notation we deduce the existence of the following fibrations. In the next result we will use the description of simple roots given in \cite[Tables]{bou}

\begin{coro}
\label{theo:U}
\label{U-fib}
\label{theo:U-non-min}
\label{U-fib-non-min}
We have the following fibrations.
\begin{enumerate}
\item If $G$ is of type $A_n$, if $\Sigma(P)=\{\a_i\}$ with $2\leq i\leq n$ and if $\Sigma(Q)=\{\iota(\a_{i-1})\}$, then $Y$ is isomorphic to $\p^{n-i+1}$ and $\vp:U\to Y$ is the vector bundle associated to the locally free sheaf $(T_{\p^{n-i+1}}(-1))^{\oplus i-1}\simeq(\Omega_{\p^{n-i+1}}^{n-i}(n-i+1))^{\oplus i-1}$.
\item If $G$ is of type $B_n$, if $\Sigma(P)=\{\a_1\}$ and if $\Sigma(Q)=\{\a_{n}\}$, then $Y$ is isomorphic to $\p^{n-1}$ and $\vp:U\to Y$ is a sequence of affine fibrations whose direction vector bundles are associated to the locally free sheaves $\co_{\p^{n-1}}(1)$ and $\Omega_{\p^{n-1}}^1(2)$.
\item If $G$ is of type $B_n$, if $\Sigma(P)=\{\a_2\}$ and if $\Sigma(Q)=\{\a_{1}\}$, then $Y$ is isomorphic to $\Q_{2n-3}$ and $\vp:U\to Y$ is a vector bundle associated to the locally free sheaf $(\co_{\Q_{2n-3}}(1))^\perp$ \emph{i.e.} to the restriction of the tautological quotient bundle on $\p^{2n-2}$ to $\Q_{2n-3}$.
\item If $G$ is of type $D_n$, if $\Sigma(P)=\{\a_1\}$ and if $\Sigma(Q)=\{\a_{n}\}$, then $Y$ is isomorphic to $\p^{n-2}$ and $\vp:U\to Y$ is the vector bundle associated to the locally free sheaf $\Omega_{\p^{n-2}}^1(2)$.
\item If $G$ is of type $D_n$, if $\Sigma(P)=\{\a_2\}$ and if $\Sigma(Q)=\{\a_{1}\}$, then $Y$ is isomorphic to $\Q_{2n-4}$ and $\vp:U\to Y$ is a vector bundle associated to the locally free sheaf $(\co_{\Q_{2n-4}}(1))^\perp$ \emph{i.e.} to the restriction of the tautological quotient bundle on $\p^{2n-3}$ to $\Q_{2n-4}$.
\item  If $G$ is of type $D_n$, if $\Sigma(P)=\{\a_n\}$ and if $\Sigma(Q)=\{\iota(\a_{n-1})\}$, then $Y$ is isomorphic to $\p^{n-1}$ and $\vp:U\to Y$ is the vector bundle associated to the locally free sheaf $\Omega_{\p^{n-1}}^{n-3}(n-2)$.
\item  If $G$ is of type $E_6$, if $\Sigma(P)=\{\a_1\}$ and if $\Sigma(Q)=\{\iota(\a_{6})\}=\{\a_1\}$, then $Y$ is isomorphic to $\Q_{8}$ and $\vp:U\to Y$ is the vector bundle associated to the locally free sheaf $E=S$ one of the two spinor bundles on $\Q_8$.
\item If $G$ is of type $E_6$, if $\Sigma(P)=\{\a_2\}$ and if $\Sigma(Q)=\{\a_1,\a_6\}$, then $Y$ is isomorphic to $\Q_{6}$ and $\vp:U\to Y$ is a sequence of affine bundles whose direction vector bundles are associated to the locally free sheaves $E=(\co_{\Q_6}(1))^\perp$ and $E'=S\oplus S'$ where $E$ is the restriction of the tautological quotient bundle of $\p^7$ to $\Q_6$ and $S$ and $S'$ are the two spinor bundles on $\Q_6$.
\item If $G$ is of type $E_7$, if $\Sigma(P)=\{\a_7\}$ and if $\Sigma(Q)=\{\a_2\}$, then $Y$ is isomorphic to $\p^{6}$ and $\vp:U\to Y$ is a sequence of affine bundles whose direction vector bundles are associated to the locally free sheaves $E=\Omega_{\p^6}^5(6)$ and $E'=\Omega_{\p^6}^2(3)$.
\item If $G$ is of type $E_7$, if $\Sigma(P)=\{\a_2\}$ and if $\Sigma(Q)=\{\a_6\}$, then $Y$ is isomorphic to $\Q_{8}$ and $\vp:U\to Y$ is a sequence of affine bundles whose direction vector bundles are associated to the locally free sheaves $E=(\co_{\Q_8}(1))^\perp$ and $E'=S^2$ where $E$ is the restriction of the tautological quotient bundle of $\p^9$ to $\Q_8$ and $S$ is one of the two spinor bundles on $\Q_8$.
\item If $G$ is of type $E_8$, if $\Sigma(P)=\{\a_8\}$ and if $\Sigma(Q)=\{\a_1\}$, then $Y$ is isomorphic to $\Q_{12}$ and $\vp:U\to Y$ is a sequence of affine bundles whose direction vector bundles are associated to the locally free sheaves $E=(\co_{\Q_{12}}(1))^\perp$ and $E'=S$ where $E$ is the restriction of the tautological quotient bundle of $\p^{13}$ to $\Q_{12}$ and $S$ is one of the two spinor bundles on $\Q_{12}$.
\item If $G$ is of type $F_4$, if $\Sigma(P)=\{\a_1\}$ and if $\Sigma(Q)=\{\a_4\}$, then $Y$ is isomorphic to $\Q_{5}$ and $\vp:U\to Y$ is a sequence of affine bundles whose direction vector bundles are associated to the locally free sheaves $E=(\co_{\Q_{5}}(1))^\perp$ and $E'=S$ where $E$ is the restriction of the tautological quotient bundle of $\p^{6}$ to $\Q_{5}$ and $S$ is the spinor bundle on $\Q_{5}$.
\end{enumerate}
\end{coro}

\begin{proo}
According to the previous Lemma, we only need to identify in each case the direction vector bundles of the affine bundles and to notice (see \cite{fourier} proof of Proposition 5) that the first affine bundle over $Y$ is always a vector bundle.

To identify the vector bundles, we only have to identify the corresponding representations. Indeed, recall from \cite[Proposition 5]{fourier} that the direction vector bundles are equivariant and correspond to the $L_{Q^{w_0}}\cap P$-representations $U^i_{Q^{w_0}}/U^{i+1}_{Q^{w_0}}\cdot (U^i_{Q^{w_0}}\cap P)$. We therefore only have to compute their highest weight which is an easy check using \cite[Tables]{bou}. Note that the representations corresponding to the sheaf $(\co_{\Q_n}(1))^\perp$ on an $n$-dimensional quadric is not irreducible as a representation of the semisimple part of $L_{Q^{w_0}}\cap P$ but is an irreducible non simple representation of $L_{Q^{w_0}}\cap P$. Note also that the above description of vector bundles for the adjoint varieties were given in \cite{rational}.
\end{proo}

%%%%%%%%%%%%%%%%%%%%%%%%%%%%%%%%%% 3 %%%%%%%%%%%%%%%%%%%%%%%%

\section{Rectriction of some homogeneous vector bundles to curves}
\label{sect-annu}

In this section we gather some known results and some new ones on the geometry of curves in the projective space
and in smooth quadrics. As explained in the previous sections, all the affine bundles $\vp:U\to Y$ we have to consider are homogeneous vector bundles over these varieties. Therefore the main aim of the section is to control the cohomology of the pull-back to curves of some equivariant vector bundles. Even if we only need these results for elliptic curves, we state them for any curve of genus $g$.

\subsection{Projective spaces}
\label{sect-secantes}

In this subsection we recall some
results on secants to a curve. The following result can be found in \cite[Proposition IV.3.8 and Theorem IV.3.9]{hartshorne}.

\begin{prop}
\label{prop-secantes}
Let $C$ be a reduced and irreducible curve in $\p^r$ such that any secant line to $C$ is a multisecant, then $C$ is a line.
\end{prop}

We will use this statement via the following corollary.

\begin{coro}
\label{coro-secantes}
Let $C$ be a reduced and irreducible curve in $\p^r$ which is not a line, then there exists a point $x$ in $C$ such that the projection from $x$ is birational from $C$ onto its image.
\end{coro}

Let $C$ be a smooth curve of genus $g$ and let $f:C\to\p^r$ be a morphism such that $f^*\co_{\p^r}(1)$ is of degree $d$. We want to prove that the pull back of the homogeneous bundles $\Omega_{\p^r}^k(k)$ through $f$ admits nice filtrations. The following lemma can be found in \cite[Remark (2) page 498]{GLP}.

\begin{lemm}
Assume that $f$ is locally injective and that $f(C)$ is non degenerate
\emph{i.e.} not contained in an hyperplane. Then, for almost all choices of $r-1$ points $(x_i)_{i\in[1,r-1]}$ in $C$, there exists a filtration
$$0=F^{r+1}\subset F^r\subset\cdots\subset F^1=f^*\Omega_{\p^r}^1(1)$$
such that $F^i/F^{i+1}=\co_C(-x_i)$ for $i\in[1,r-1]$ and $F^r/F^{r+1}=f^*\co_{\p^r}(-1)\otimes\co_C(\sum_i x_i)$.
\end{lemm}

\begin{proo}
Let us set $M=f^*\Omega_{\p^r}^1(1)$ and $L=f^*\co_{\p^r}(1)$. Let us also denote by $V^\vee$ the vector space $H^0(\p^r,\co_{\p^r}(1))$. Because $f(C)$ is non degenerate, we may consider $V^\vee$ as a subspace of $H^0(C,L)$. We have an exact sequence $0\to M\to V^\vee\otimes\co_C\to L\to0$.

We proceed by induction on $r$. Let us choose a point $x_1$ on $C$. Let us denote by $W^\vee$ the intersection $V^\vee\cap H^0(C,L(-x_1))$ in the vector space $H^0(C,L)$. If $f(x_1)$ is smooth, the map $W^\vee\otimes\co_C\to L(-x_1)$ is surjective and we have the commutative diagram:
$$\xymatrix{&&W^\vee\otimes\co_C\ar[d]\ar[r]&L(-x_1)\ar[d]&\\
0\ar[r]&M\ar[r]\ar[d]&V^\vee\otimes\co_C\ar[r]\ar[d]& L\ar[r]\ar[d]&0\\
0\ar[r]&\co_C(-x_1)\ar[r]&\co_C\ar[r]&\co_{x_1}\ar[r]&0.}$$
The snakes lemma implies that the left most vertical map is surjective
and we get the term $F^2$ of the filtration by taking the kernel of
this map which is also the kernel of the map $W^\vee\otimes\co_C\to
L(-x_1)$. Replacing $V$ by $W$ and $L$ by $L(-x_1)$ we get a morphism
$f_2$ from $C$ to $\p^{r-1}$ and $F^2=f_2^*\Omega_{\p^{r-1}}(1)$. The
morphism $f_2$ is the composition of $f$ with the projection of center
$x_1$ and because, in characteristic zero, general secants are not
multisecant (see Corollary \ref{coro-secantes}) we have that $f_2$ is also generically
injective for $x_1$ general. We may apply then the induction
hypothesis. For $r=1$, the result is true.
\end{proo}

\begin{rema}
  Note that in the above filtration, all the vector bunndles $F^i$ are
  contructed as tautological subbundles therefore their dual are
  globally generated.
\end{rema}

\begin{coro}
Assume that $f(C)$ is non degenerate, then for almost all choices of $r-1$ points $(x_i)_{i\in[1,r-1]}$ in $C$, there exists a filtration
$$0=G^{s+1}\subset G^s\subset\cdots\subset G^1=f^*\Omega_{\p^r}^k(k)$$
such that $G^i/G^{i+1}=\co_C(-(x_{i_1}+\cdots+x_{i_k}))$ with $\{i_1,\cdots,i_k\}\subset[1,r-1]$ or $G^i/G^{i+1}=f^*\co_{\p^r}(-1)\otimes\co_C(x_{i_1}+\cdots+x_{i_{r-k}})$ with $\{i_1,\cdots,i_{r-k}\}\subset[1,r-1]$.
\end{coro}

\begin{proo}
We simply take the $k$-th exterior product of the previous filtration.
\end{proo}

\begin{coro}
\label{coro-van}
Assume that $f(C)$ is non degenerate. Let $g$ be the genus of $C$ and $d$ be the degree of $f^*\co_{\p^r}(1)$. If $g\leq r-k$, then we have the vanishing $H^1(C,f^*\Omega_{\p^r}^k(k+1))=0$.
\end{coro}

\begin{proo}
Indeed, for $g\leq r-k$ and for general points $(x_i)_{i\in[1,r-1]}$ on $C$, we have the vanishing $H^1(C,\co_C(x_{i_1}+\cdots+x_{i_{r-k}}))=0$. Note that because $f(C)$ is non degenerate, we must have the inequality $g\leq d-r\leq d-k$ and this condition implies the following vanishing $H^1(C,f^*\co_{\p^r}(1)\otimes\co_C(-(x_{i_1}+\cdots+x_{i_{k}})))=0$.
\end{proo}

Let us assume that $C$ is elliptic. Then even if $f(C)$ is degenerate we may obtain results. Indeed, if $f(C)$ is degenerate, there is a minimal linear space $H$ in $\p^r$ of codimension say $a$ containing $f(C)$ and such that $f(C)$ is non degenerate in $H$. Furthermore we have the exact sequence $0\to\co_C^a\to f^*\Omega_{\p^r}^1(1)\to f^*\Omega_H^1(1)\to 0$ from which we deduce a filtration of $f^*\Omega_{\p^r}^1(1)$. In particular we get the following result.

\begin{prop}
\label{prop-proj-ell}
Assume that $C$ is elliptic.
If $f(C)$ is non degenerate in a linear subspace of codimension
$a$, then we have the equality $\dim
H^1(C,f^*\Omega_{\p^r}^k(k+1))=\binom{a}{r-k}$ (with
$\binom{a}{r-k}=0$ if $a<r-k$).

In particular, if we have the inequality $\dim
H^1(C,f^*\Omega_{\p^r}^k(k+1))>0$, then there exists an integer $a\geq
r-k$
such that $\dim H^1(C,f^*\Omega_{\p^r}^k(k+1))=\binom{a}{r-k}$ and
$f(C)$ in non degenerate in a linear subspace of codimension $a$.
\end{prop}

\begin{proo}
Taking the $k$-th exterior power of the exact sequence
$0\to\co_C^a\to f^*\Omega_{\p^r}^1(1)\to f^*\Omega_H^1(1)\to 0$ and
tensoring it by $f^*\co_{\p^r}(1)$, we get a filtration of
$f^*\Omega_{\p^r}^k(k+1)$ by the vector bundles
$\Lambda^u(\co_C^a)\otimes f^*\Omega^v_{H}(v+1)$ for $u+v=k$ and $u,v$
non negative. By Corollary \ref{coro-van}, the first cohomology groups
of these bundles vanish except maybe for $v=0$ or $v=r-a$. For $v=0$ we have
$\Lambda^u(\co_C^a)\otimes
f^*\Omega^v_{H}(v+1)=\Lambda^k(\co_C^a)\otimes f^*\co_{\p^r}(1)$
and for $v=r-a$ we have $\Lambda^u(\co_C^a)\otimes
f^*\Omega^v_{H}(v+1)=\Lambda^{k+a-r}(\co_C^a)$. The first cohomology
group in the first case vanishes while the first cohomology group in the
second case has dimension $\binom{a}{r-k}$. We thus have an inequality
$\dim H^1(C,f^*\Omega^k_{\p^r}(k+1))\leq\binom{a}{r-k}$.

For $a<r-k$, we are done. If $a\geq r-k$, let us consider the
identification
$\Omega_{\p^r}^k(k+1)=(\Omega_{\p^r}^{r-k}(r-k))^\vee$. By the above
exact sequence, we have a trivial subbundle of rank $a$ of
$f^*\Omega_{\p^r}^1(1)$ therefore we have a trivial subbundle of rank
$\binom{a}{r-k}$ in $f^*\Omega^{r-k}_{\p^r}(r-k)$ (take the $r-k$
exterior power). By duality, there is a trivial quotient bundle of
rank $\binom{a}{r-k}$ of $f^*\Omega_{\p^r}^k(k+1)$.

\begin{lemm}
\label{lemm-ell-sections}
  Let $E$ be a globally generated vector bundle on $C$ such that $E$
  has a trivial quotient of rank $n$, then $\dim H^1(C,E)\geq n$.
\end{lemm}

\begin{proo}
  By Serre duality we have an isomorphism $H^1(C,E)\simeq
  H^0(C,E^\vee)^\vee$. We therefore need to prove that $\dim
  H^0(C,E^\vee)\geq n$. But we have a surjection $E\to\co_C^n$ giving
  rise to an injection $\co_C^n\to E^\vee$ and the result follows.
\end{proo}

We conclude
using the previous lemma and because $\Omega_{\p^r}^k(k+1)$ is globally
generated.
\end{proo}

\subsection{Spinor bundles on quadrics}

In this section we consider the spinor bundles on quadrics and their restriction to elliptic curves.

\begin{prop}
Let $\Q_n$ be a smooth $n$-dimensional quadric and let $E$ be a spinor bundle on $\Q_n$. Let $f:C\to \Q_n$ be a morphism from a smooth irreducible curve to $\Q_n$. Then the pull-back $f^*E$ is isomorphic to a direct sum ${\co_C}^{2^\beta}\oplus F$ where $F$ has no trivial factor if and only if $f$ factors through an isotropic subspace of (linear) dimension $[n/2]+1-\a$ with
$$\a=\left\{\begin{array}{ll}
\beta & \textrm{for $n$ odd},\\
\beta+1 & \textrm{for $n$ even and $\beta>0$},\\
0 & \textrm{for $n$ even and $\beta=0$}\\
\end{array}\right.$$
and in the last case, the maximal isotropic subspace has to be an element of a fixed component of the two connected components of maximal isotropic subspaces.
\end{prop}

\begin{proo}
We shall here use the results of G. Ottaviani \cite{ottaviani} on the spinor bundles. Let us write $h=\lceil n/2\rceil-1$ and $p=[n/2]+1$. There is an embedding $\varphi:\Q_n\to\G(2^h,2^{h+1})$ such that $E=\varphi^*K$ with $K$ the tautological subbundle in the Grassmann variety $\G(2^h,2^{h+1})$ of $2^h$ dimensional subspaces in $\C^{2^{h+1}}$.

\begin{lemm}
Let $x$ and $y$ be points on $\Q_n$, then we have the equality $\C^{2^{h+1}}=K_{\varphi(x)}\oplus K_{\varphi(y)}$ if and only if $x$ and $y$ are not orthogonal.
\end{lemm}

\begin{proo}
We shall consider the dual assertion: the map $\C^{2^{h+1}}\to K_{\varphi(x)}^\vee\oplus K_{\varphi(y)}^\vee$ is an isomorphism if and only if $x$ and $y$ are not orthogonal. To test this assertion, we can restrict the bundle $K$ to some subvariety containing both $x$ and $y$. If $x$ and $y$ are not orthogonal, then we may choose a smooth conic $c:\p^1\to Q_n\to \G(2^h,2^{h+1})$ passing through $x$ and $y$. By \cite[Theorem 1.4]{ottaviani}, the pull back $c^*K^\vee$ is isomorphic to $\co_{\p^1}(1)^{2^h}$. The surjectivity condition holds.

Conversely, if $x$ and $y$ are orthogonal, then there exists a maximal isotropic subspace $V_p$ (if $n$ is even we can choose $V_p$ to be in any of the two conneced components of maximal isotropic subspaces) such that $x$ and $y$ are contained in $\p(V_p)$. By \cite[Theorems 2.5 and 2.6]{ottaviani} the restriction of $K$ to $\p(V_p)$ contains a trivial factor (if $n$ is even this is true only for $V_p$ in one of the two connected components). Therefore the surjectivity condition does not hold.
\end{proo}

Let us return to the proof of the proposition. Remark that $f^*E$ has a trivial factor if and only if for all $x\in C$ the subspaces $K_{\varphi(x)}$ have a common intersection in $\C^{2^{h+1}}$. This occurs only if for all $x$ and $y$ in $C$, the points $f(x)$ and $f(y)$ are orthogonal. In particular the linear subspace spanned by $f(C)$ has to be isotropic.  Therefore if the curve $f(C)$ is not contained in an isotropic subspaces then $f^*E$ has no trivial factor.

Let us assume that $f^*E$ has a trivial factor. Then $f$ factors through an isotropic subspace $V_p$. Furthermore by \cite[Theorems 2.5 and 2.6]{ottaviani}, we know that the restriction $E_{\p(V_p)}$ is isomorphic to one of the following direct sums (the first case occurs for $n$ odd while the other two cases occur for $n$ even):
$$\bigoplus_{k=0}^{p}\Omega_{\p(V_p)}^k(k),\ \ \ \bigoplus_{k=0,\ k\textrm { odd}}^{p}\Omega_{\p(V_p)}^k(k)\ \ \textrm{ or }\ \ \bigoplus_{k=0,\ k\textrm { even}}^{p}\Omega_{\p(V_p)}^k(k).$$
By the previous results on the projective space, we conclude that if $f(C)$ factors through a linear subspace of linear dimension $p-\a$  and if there is a unique such factor then $V_p$ is in the connected component of $\G_Q(p,2p)$ corresponding to the last decomposition above. The number of trivial factors in $f^*E$ is then $2^\beta$ with $\beta$ and $\a$ as above.
\end{proo}

For elliptic curves, applying the same techniques as in the case of the projective space, we obtain the following result.

\begin{coro}
\label{spinor}
Let $\Q_n$ be a smooth $n$-dimensional quadric and let $E$ be a spinor bundle on $\Q_n$. Let $f:C\to \Q_n$ be a morphism from a smooth irreducible elliptic curve to $\Q_n$. Then we have the equality $\dim H^1(C,f^*E)=2^\beta$ if and only if $f$ factors through an isotropic subspace of (linear) dimension $[n/2]+1-\a$ with
$$\a=\left\{\begin{array}{ll}
\beta & \textrm{for $n$ odd},\\
\beta+1 & \textrm{for $n$ even and $\beta>0$},\\
0 & \textrm{for $n$ even and $\beta=0$}\\
\end{array}\right.$$
and in the last case, the maximal isotropic subspace has to be an element of a fixed component of the two connected components of maximal isotropic subspaces.
\end{coro}

%%%%%%%%%%%%%%%%%%%%%%%%%%%%%%%%%% 4 %%%%%%%%%%%%%%%%%%%%%%%%

\section{Irreducibility}
\label{sect-irr}

We are now in position to prove our main result. Let $C$ be a smooth elliptic curve. We proceed via a case-by-case analysis to finish the proof of Theorem \ref{main}.

\subsection{Projective space}

Our proof is based on the fact that the irreducibility of the scheme of morphisms to a projective space is well known. Let us first recall this fact. Let $X$ be the projective space $\p(V)$ with $V$ of dimension $n$.

\begin{fact}
Let $d\geq2$ be an integer, the scheme $\MorC{d}{X}$ is irreducible of dimension $dn$.
\end{fact}

\begin{proo}
The scheme of morphism is described by simple Brill-Noether data: a line bundle ${\cal L}$ on $C$ of degree $d$ and a linear map from $V^\vee$ to $H^0(C,{\cal L})$. But for $d$ positive, we have $H^1(C,{\cal L})=0$ therefore $H^0(C,{\cal L})$ is of constant dimension $d$ for all ${\cal L}$ and the result follows from the fact that ${\rm Jac}^d(C)$, the Jacobian of degree $d$ line bundles on $C$, is irreducible: the variety $\MorC{d}{\p(V)}$ is a locally trivial bundle over ${\rm Jac}^d(C)$ with fiber $\p(V\otimes H^0(C,{\cal L}))$ over ${\cal L}$.
\end{proo}

\subsection{Grassmann varieties}

We now deal with Grassmann varieties. Note that in this case the result was proved by A. Brugui{\`e}res \cite{bruguieres}. We include here a proof that will serve as a model for the other cases. Let $X$ be the Grassmann variety $\G(p,n)$ of $p$-dimensional vector spaces in a vector space of dimension $n$.

\begin{prop}
The scheme $\MorC{d}{X}$ is irreducible of dimension $nd$ for $d\geq n$.
\end{prop}

\begin{proo}
Consider the open subset $U$ given by Theorem \ref{U-fib}. We have a morphism $\phi:U\to Y$ with $Y\simeq\p^{n-p}$ which is the vector bundle associated to the locally free sheaf
$$E=T_{\p^{n-p}}(-1)\otimes \co_{\p^{n-p}}^{p-1}=(\Omega_{\p^{n-p}}^{n-p-1}(n-p))^{p-1}.$$
By Proposition \ref{ouvert}, we only need to prove the irreducibility of the scheme $\MorC{d}{U}$. For this we study the map $\phi$. Let us denote by ${\bf H}_a$ the locally closed subset of $\MorC{d}{Y}$ of maps whose image is contained in a linear subspace of codimension $a$ and not in a smaller linear subspace. The dimension of this locally closed subset is $a(n-p+1-a)+d(n-p+1-a)$. Indeed, there is a morphism ${\bf H}_a\to \G(n-p+1-a,n-p+1)$ sending to a morphism the unique linear subspace of codimension $a$ containing its image. The fibers are isomorphic to $\MorC{d}{\p^{n-p-a}}$.

Now the inverse image under the map $\Phi:\MorC{d}{U}\to\MorC{d}{Y}$ induced by $\phi$ of ${\bf H}_a$ is irreducible of dimension $\dim {\bf H}_a+\dim H^0(C,f^*E)$ for $f\in{\bf H}_a$. By Proposition \ref{prop-proj-ell} we have $\dim H^0(C,f^*E)=(d\deg(E)+a)(p-1)=d(p-1)+a(p-1)$. In particular for $d\geq n$, the dimension of $\Phi^{-1}({\bf H}_a)$ is $(d+a)(n-a)=dn-a(d-n)-a^2<dn$ for $a>0$. The locally closed subsets $\Phi^{-1}({\bf H}_a)$ for $a>0$ are therefore of dimension smaller than the dimensions (at least $dn$) of irreducible components of $\MorC{d}{U}$ thus the irreducible components of $\Mor{d}{U}$ are those of $\Phi^{-1}({\bf H}_0)$ which is irreducible of dimension $dn$.
\end{proo}

\begin{rema}
For smaller degrees, the scheme $\MorC{d}{X}$ is not irreducible (except in degree 2). Its irreducible components were described in \cite{bruguieres}.
\end{rema}

\subsection{Quadrics}

We shall now deal with smooth quadrics. Note that in this case the irreducibility of the scheme of morphisms was proved by see E. Ballico \cite{ballico} for quadrics of dimension $n\geq 6$ for curves of genus $g$ with $g\leq\lfloor\frac{n-1}{2}\rfloor$ and for the degree $d\geq 2g-1$. Let $X$ be a smooth quadric of dimension $n$.

\begin{prop}\label{prop:quad}
Assume that $\dim X\neq3$, then the scheme $\MorC{d}{X}$ is irreducible of dimension $nd$ for $d\geq 3$.
\end{prop}

\begin{proo}
We proceed as for the Grassmann variety and consider the open subset $U$ given by Theorem \ref{U-fib}. We have a sequence of affine bundles $\phi:U\to Y$ where $Y\simeq\p^r$ with $r=\lfloor n/2\rfloor$ and the direction vector bundles are associated to the locally free sheaves $E=\Omega^1_{\p^r}(2)$ for $n$ even (here $\vp$ is a vector bundle) and to $E=\Omega^1_{\p^r}(2)$ and $E'={\cal O}_{\p^r}(1)$ for $n$ odd.

By Proposition \ref{prop-proj-ell} the locally free sheaf $E$ satisfies $H^1(C,f^*E)=0$ for any map $f:C\to Y$ of degree $d\geq3$ as soon as $f$ does not factor through a line. In this latter case we have $\dim H^1(C,f^*E)=1$ and we conclude by Proposition \ref{prop-vb}.
\end{proo}

\begin{rema}
(\i) Note that if $\dim X\leq2$, then $X$ is a product of projective lines and the result follows from that case. However, if $\dim X=3$, then in the above proof we get $r=1$ and $\Omega^1_{\p^r}(2)=\co_{\p^r}$ and we do not get the vanishing condition.

(\i\i) Note that for $d=2$ the scheme $\MorC{d}{X}$ is also irreducible but not of the expected dimension since any degree 2 morphism factors through a line. The dimension is $2n+1$.
\end{rema}

Even if we cannot prove, with our method, Proposition \ref{prop:quad} when $X$ is 3-dimensional, the result is still true.

\begin{prop}\label{prop:q3}
Let $X=\Q_3$ be a smooth quadric of dimension~3, then the scheme ${\bf Hom}_{d}(C,X)$ is irreducible of dimension $3d$ for $d\geq 3$.
\end{prop}

We give a proof of this result for $d\geq 4$ in the appendix, using a different method. The case where $d=3$ is given by the following remark.

\begin{rema}
Any morphism of degree 2 or 3 from an elliptic curve to $\Q_3$ has to factor through a line. This is clear for degree 2, and for degree 3 the morphism factors through a plane. But as $\Q_3$ contains no plane, the image is contained in the intersection of $\Q_3$ with a plane. It has to be a line for degree reasons. In particular the scheme of morphisms ${\bf Hom}_{d}(C,X)$ is irreducible of dimension $3+2d$ for $d=2$ or $d=3$.
\end{rema}

\subsection{Maximal orthogonal Grassmann varieties}

Let $X$ be a maximal orthogonal Grassmann variety $\G_Q(n,2n)$. The following result already appeared in the unpublished work of the second author \cite{arxiv} with a different method.

\begin{prop}
The scheme $\MorC{d}{X}$ is irreducible of dimension $2(n-1)d$ for $d\geq n-1$.
\end{prop}

\begin{proo}
The map $\phi:U\to Y$ is the vector bundle over $Y\simeq \p^{n-1}$ associated to the locally free sheaf $E=\Lambda^2(T_{\p^{n-1}}(-1))=\Omega^{n-3}_{\p^{n-1}}(n-2)$. By Proposition \ref{prop-proj-ell}, the fibers of $\Phi:\MorC{d}{U}\to \MorC{d}{Y}$ over the locally closed subset ${\bf H}_a$ of morphisms whose image is non-degenerate in a linear subspace of codimension $a$ are of dimension $d(n-2)+\binom{a}{2}$. Therefore we have $\dim p^{-1}({\bf H}_a)=2d(n-1)+a(n-a)+\binom{a}{2}-ad$. This dimension is strictly smaller than the expected dimension $2d(n-1)$ for $a>0$ and $d\geq n-1$. The result follows.
\end{proo}

\begin{coro}
\label{coro-gq}
For $d\in[2,n]$, the scheme $\MorC{d}{X}$ is irreducible of dimension
$$2(n-1)d+\frac{1}{2}(n-d)(n-d-1).$$
\end{coro}

\begin{proo}
Let $f:C\to X$ be a morphism of degree $d$, then the kernel $\ker(f)$ of $f$ \emph{i.e.} the intersection of all subspaces $f(x)$ for $x\in C$ is of dimension at least $n-d$ (see \cite{buch} for other applications of this definition). Indeed, consider the tautological exact sequence
$$0\to K\to {\cal O}_X^{2n}\to Q\to 0$$
on $X$ where $K$ and $Q$ are the tautological subbundle and quotient bundle. On $X$ we have the identification $Q\simeq K^\vee$. Pulling back this sequence to $C$ via $f$ we get an exact sequence on the curve. Taking cohomology we obtain:
$$0\to H^0(C,K)\to \C^{2n}\to H^1(C,K)\to H^1(C,K)\to \C^{2n}\to H^0(C,K)\to 0.$$
Let us denote by $h^0$ and $h^1$ the dimensions of $H^0(C,K)$ and $H^1(C,K)$. Since the map $f$ is of degree $d$ in $X$, the degree of $K$ is $-2d$ and we have the equality $h^0=h^1-2d$. We deduce the inequality $h^0\geq n-d$. This implies that $\ker(f)$ is of dimension at least $n-d$. Note also that for a general curve the kenel is of dimension exactly $n-d$.

Let us consider the following incidence variety
$$I=\{(f,V_{n-d})\in \MorC{d}{X}\times\G_Q(n-d,2n)\ /\ \textrm{ for all $x\in C$ we have $V_{n-d}\subset f(x)\subset V_{n-d}^\perp$}\}.$$
Here $\G_Q(n-d,2n)$ is the Grassmann variety of isotropic subspaces of $\C^{2n}$ of dimension $n-d$. The above condition on $f$ simply translates in $V_{n-d}\subset \ker(f)$. Note that the second inclusion is implied by the first one since we have $f(x)=f(x)^\perp$. The projection $p:I\to \MorC{d}{X}$ is surjective and birational (because for a general curve the kernel is of dimension $n-d$) while the map $q:I\to \G_Q(n-d,2n)$ is an equivariant locally trivial fibration with fibers $\MorC{d}{Y}$ with $Y\simeq\G_Q(d,2d)$. By the previous proposition these fibers are irreducible of dimension $2d(d-1)$. The result follows.
\end{proo}

\subsection{Adjoint varieties in type $B$ and $D$}

Let us consider the variety $X=\G_Q(2,n)$ of isotropic planes in a vector space of dimension $n$ endowed with a non-degenerate bilinear form.

\begin{prop}
The scheme $\MorC{d}{X}$ is irreducible of dimension $(n-3)d$ for $d\geq n-1$.
\end{prop}

\begin{proo}
The map $\phi:U\to Y$ is a vector bundle over $Y\simeq \Q_{n-4}$ associated to the locally free sheaf $E=(\co_{\Q_{n-4}}(1))^\perp$ where $\co_{\Q_{n-4}}(-1)$ is the tautological subbundle on  $\Q_{n-4}$.

Here we have to compute the dimension of the scheme of morphisms from an elliptic curve to a cone (the intersection of $\Q_{n-4}$ with a linear subspace of $\p^{n-3}$ in our case). Let $\mathfrak{C}$ be a cone in $\p^{k+r-1}$ associated to a quadratic form of rank $r$ (and therefore of kernel $K$ of dimension $k$).

\begin{lemm}
Consider the open subset in the scheme of degree $d$ morphisms from an elliptic curve $C$ to $\mathfrak{C}$ of morphisms whose image is non-degenerate. Then this open subset is of dimension at most $d(k+r-2)$.
\end{lemm}

Remark that, if $k=0$, the cone $\mathfrak{C}$ is a smooth quadric of dimension $(k+r-2)$ and we already know the result with equality.

\begin{proo}
Consider the incidence variety $I=\{(V_1,V_{k+1})\in\mathfrak{C}\times\G_Q(k+1,k+r)\ /\ V_1\subset V_{k+1}\supset K\}$. There is a diagram
$$\xymatrix{I\ar[r]^p\ar[d]_q & \Q_{r-2}\\
\mathfrak{C}. & \\}$$
The map $q$ is simply the blow-up of the vertex of the cone. If $f:C\to\mathfrak{C}$ has a non degenerate image then its image is not contained in the vertex $\p(K)$ of the cone and we may consider the composition of $f$ with the projection from $K$. This defines a morphism $f':C\to\Q_{r-2}$ which can be obtained by lifting $f$ to $f_I:C\to I$ and composing with $p$. The degree of this morphism is $d-x$ where $x$ is the multiplicity of the intersection of $f(C)$ with the vertex of the cone.

To compute the dimension of our open set, we consider the map defined by $f\mapsto f'$. Its image is contained in $\MorC{d-x}{\Q_{r-2}}$ which is of dimension $(d-x)(r-2)$. The fiber is given by the morphisms ${\cal L}\to {f'}^*V_{k+1}$ where ${\cal L}$ is of degree $d$ and $V_{k+1}$ is the tautological subbundle in $\Q_{r-2}\simeq\{V_{k+1}\in\G_Q(k+1,k+r)\ /\ K\subset V_{k+1}\}$. Because the projection by $K$ gives a curve of degree $d-x>0$ we have an extension $0\to K\otimes\co_C\to {f'}^*V_{k+1}\to {\cal L}'\to 0$ with ${\cal L}'$ of degree $d-x$. This extension has to split, thus the fibres of the map are given by $\p{\rm Hom}({\cal L},K\otimes\co_C\oplus{\cal L}'))$ which is of dimension $kd+x$. The dimension is therefore at most $d(k+r-2)-x(r-3)$ proving the result for $r\geq 3$. But if $r\leq2$, then the cone is a union of hyperplane thus the curve is degenerate.
\end{proo}

By Proposition \ref{prop-proj-ell}, since $E$ is the restriction to $\Q_{n-4}$ of $\Omega_{\p^{n-3}}^{n-4}(n-3)$, the fibers of $\Phi:\MorC{d}{U}\to \MorC{d}{Y}$ over the locally closed subset ${\bf H}_a$ of morphisms whose image is non-degenerate in a linear subspace of codimension $a$ are of dimension $d+a$. Therefore the dimension of $\Phi^{-1}({\bf H}_a)$ is at most the sum of $d+a$, the dimension of the Grassmann variety $\G(n-2-a,n-2)$ and the dimension of the scheme of degree $d$ morphisms from an elliptic curve $C$ to a cone $\mathfrak{C}$ in $\p^{n-3-a}$, so that $\dim \Phi^{-1}({\bf H}_a)\leq d(n-a-3)+a(n-2-a)+a$. This dimension is strictly smaller than the expected dimension $(n-4)d$ for $a>0$ and $d\geq n-1$. The result follows.
\end{proo}

\begin{rema}
The above argument shows that for $d=n-2$, the scheme of morphisms has 2 irreducible components of the expected dimension.
\end{rema}

\subsection{Type $E_6$}

\subsubsection{The Cayley plane}

Let $X$ be isomorphic to $E_6/P_1$.

\begin{prop}
The scheme $\MorC{d}{X}$ is irreducible of dimension $12d$ for $d\geq 3$.
\end{prop}

\begin{proo}
The map $\phi:U\to Y$ is the vector bundle over $Y\simeq \Q_8$ associated to the locally free sheaf $E=S$ with $S$ one of the two spinor bundles on $\Q_8$. By Corollary \ref{spinor}, the fibers of $p:\MorC{d}{U}\to \MorC{d}{Y}$ over the locally closed subset ${\bf H}_a$ of morphisms whose image is non degenerate in an isotropic linear subspace of codimension $a$ are of dimension $4d+2^{\max(0,a-6)}$ with $a\geq 5$. Therefore we have the equality $\dim \Phi^{-1}({\bf H}_a)=12d+a(10-a)-\frac{(10-a)(11-a)}{2}+2^{\max(0,a-6)}-ad$. This dimension is strictly smaller than the expected dimension $12d$ for $a>0$ and $d\geq 3$. The result follows.
\end{proo}

\subsubsection{The adjoint variety}

Let $X$ be isomorphic to $E_6/P_2$.

\begin{prop}\label{prop:E6/P2}
The scheme $\MorC{d}{X}$ is irreducible of dimension $11d$ for $d\geq 9$.
\end{prop}

\begin{proo}
The map $\phi:U\to Y$ is a sequence of affine bundles over $Y\simeq \Q_{6}$ associated to the locally free sheaves $E=(\co_{\Q_6}(1))^\perp$ (where $\co_{\Q_6}(-1)$ is the tautological subbundle on $\Q_6$) and $E'=S\oplus S'$ where $S$ and $S'$ are the two spinor bundles on $\Q_6$. Note that the bundle $E$ is the restriction of the tautological quotient bundle of $\p^7$ to $\Q_6$. By Proposition \ref{prop-proj-ell} and Corollary \ref{spinor}, the fibers of $\Phi:\MorC{d}{U}\to \MorC{d}{Y}$ over the locally closed subset ${\bf H}_a$ of morphisms whose image is non degenerate in a linear subspace of codimension $a$ are of dimension at most $5d+a$ for $a\leq 3$ or $a\geq 4$ and the subspace containing the curve is non isotropic and $5d+a+2^{\max(0,a-4)}$ for $a\geq 4$ and the subspace containing the curve is isotropic. Therefore we have the equality $\dim \Phi^{-1}({\bf H}_a)=d(11-a)+a(8-a)+a$ for $a\leq 3$ or $a\geq 4$ and the subspace containing the curve is non isotropic and $d(13-a)+a(8-a)-\frac{(8-a)(9-a)}{2}+a+2^{\max(0,a-4)}$ for $a\geq 4$ and the subspace containing the curve is isotropic. This dimension is strictly smaller than the expected dimension $11d$ for $a>0$ and $d\geq 9$. The result follows.
\end{proo}

\subsection{Type $E_7$}

\subsubsection{The Freudenthal variety}

Let $X$ be isomorphic to $E_7/P_7$.

\begin{prop}\label{prop:E7/P7}
The scheme $\MorC{d}{X}$ is irreducible of dimension $18d$ for $d\geq 8$.
\end{prop}

\begin{proo}
The map $\phi:U\to Y$ is a sequence of affine bundles associated to the locally free sheaves $E=\Omega^5_{\p^6}(6)$ and $E'=\Omega^2_{\p^6}(3)$. By Proposition \ref{prop-proj-ell}, the fibers of $\Phi:\MorC{d}{U}\to \MorC{d}{Y}$ over the locally closed subset ${\bf H}_a$ of morphisms whose image is non degenerate in a linear subspace of codimension $a$ are of dimension at most $11d+a+\binom{a}{4}$. Therefore we have $\dim \Phi^{-1}({\bf H}_a)=18d+a(7-a)+a+\binom{a}{4}-ad$. This dimension is strictly smaller than the expected dimension $18d$ for $a>0$ and $d\geq 8$. The result follows.
\end{proo}

\subsubsection{The adjoint variety}

Let $X$ be isomorphic to $E_7/P_1$.

\begin{prop}\label{prop:E7/P1}
The scheme $\MorC{d}{X}$ is irreducible of dimension $17d$ for $d\geq 11$.
\end{prop}

\begin{proo}
The map $\phi:U\to Y$ is a sequence of affine bundles over $Y\simeq \Q_8$ associated to the locally free sheaves $E=(\co_{\Q_8}(1))^\perp$ where $\co_{\Q_8}(-1)$ is the tautological subbundle on $\Q_8$ and $E'=S^2$ where $S$ is one of the two spinor bundles on $\Q_8$.
 By Proposition \ref{prop-proj-ell} and Corollary \ref{spinor}, the fibers of $\Phi:\MorC{d}{U}\to \MorC{d}{Y}$ over the locally closed subset ${\bf H}_a$ of morphisms whose image is non degenerate in a linear subspace of codimension $a$ are of dimension at most $9d+a$ for $a\leq 4$ or $a\geq 5$ and the subspace containing the curve is non isotropic and $9d+a+2^{\max(0,a-6)}$ for $a\geq 5$ and the subspace containing the curve is isotropic. Therefore we have the equality $\dim \Phi^{-1}({\bf H}_a)=d(17-a)+a(10-a)+a$ for $a\leq 4$ or $a\geq 5$ and the subspace containing the curve is non isotropic and $d(19-a)+a(10-a)-\frac{(10-a)(11-a)}{2}+a+2^{\max(0,a-6)}$ for $a\geq 5$ and the subspace containing the curve is isotropic. This dimension is strictly smaller than the expected dimension $17d$ for $a>0$ and $d\geq 11$. The result follows.
\end{proo}

\subsection{Adjoint variety of type $E_8$}

Let $X$ be isomorphic to $E_8/P_8$.

\begin{prop}\label{prop:E8/P8}
The scheme $\MorC{d}{X}$ is irreducible of dimension $29d$ for $d\geq 15$.
\end{prop}

\begin{proo}
The map $\phi:U\to Y$ is a sequence of affine bundles over $Y\simeq \Q_{12}$ associated to the locally free sheaves $E=(\co_{\Q_{12}}(1))^\perp$ where $\co_{\Q_{12}}(-1)$ is the tautological subbundle on $\Q_{12}$ and $E'=S$ where $S$ is one of the two spinor bundles on $\Q_{12}$.
 By Proposition \ref{prop-proj-ell} and Corollary \ref{spinor}, the fibers of $\Phi:\MorC{d}{U}\to \MorC{d}{Y}$ over the locally closed subset ${\bf H}_a$ of morphisms whose image is non degenerate in a linear subspace of codimension $a$ are of dimension at most $17d+a$ for $a\leq 6$ or $a\geq 7$ and the subspace containing the curve is non isotropic and $17d+a+2^{\max(0,a-8)}$ for $a\geq 7$ and the subspace containing the curve is isotropic. Therefore we have the equality $\dim \Phi^{-1}({\bf H}_a)=d(29-a)+a(14-a)+a$ for $a\leq 6$ or $a\geq 7$ and the subspace containing the curve is non isotropic and $d(31-a)+a(14-a)-\frac{(14-a)(15-a)}{2}+a+2^{\max(0,a-8)}$ for $a\geq 7$ and the subspace containing the curve is isotropic. This dimension is strictly smaller than the expected dimension $29d$ for $a>0$ and $d\geq 15$. The result follows.
\end{proo}

\subsection{Adjoint variety of type $F_4$}

Let $X$ be isomorphic to $F_4/P_1$.

\begin{prop}\label{prop:F4/P1}
The scheme $\MorC{d}{X}$ is irreducible of dimension $11d$ for $d\geq 8$.
\end{prop}

\begin{proo}
The map $\phi:U\to Y$ is a sequence of affine bundles over $Y\simeq \Q_{5}$ associated to the locally free sheaves $E=(\co_{\Q_{5}}(1))^\perp$ where $\co_{\Q_{5}}(-1)$ is the tautological subbundle on $\Q_{5}$ and $E'=S$ where $S$ is the spinor bundle on $\Q_{5}$.
By Proposition \ref{prop-proj-ell} and Corollary \ref{spinor}, the fibers of $\Phi:\MorC{d}{U}\to \MorC{d}{Y}$ over the locally closed subset ${\bf H}_a$ of morphisms whose image is non degenerate in a linear subspace of codimension $a$ are of dimension at most $6d+a$ for $a\leq 3$ or $a\geq 4$ and the subspace containing the curve is non isotropic and $6d+a+2^{\max(0,a-5)}$ for $a\geq 4$ and the subspace containing the curve is isotropic. Therefore we have the equality $\dim \Phi^{-1}({\bf H}_a)=d(11-a)+a(7-a)+a$ for $a\leq 3$ or $a\geq 4$ and the subspace containing the curve is non isotropic and $d(13-a)+a(7-a)-\frac{(7-a)(8-a)}{2}+a+2^{\max(0,a-5)}$ for $a\geq 4$ and the subspace containing the curve is isotropic. This dimension is strictly smaller than the expected dimension $1d$ for $a>0$ and $d\geq 8$. The result follows.
\end{proo}

\begin{rema}
Note that, in Propositions \ref{prop:E6/P2}, \ref{prop:E7/P7}, \ref{prop:E7/P1}, \ref{prop:E8/P8} and \ref{prop:F4/P1}, it is possible that one can improve the bound on the dimension of $\Phi^{-1}({\bf H}_a)$. Indeed if the affine bundle restricted to curve in ${\bf H}_a$ is not a vector bundle then this bound will be strictly smaller. We therefore expect a better bound $d(X)$ in these situations.
\end{rema}

\section{Stable maps}\label{sect-stable}

Let us consider the moduli space $\MK{1,d}{X}$ of stable maps with source an elliptic curve. Recall the following general result obtained in \cite{KP}.

\begin{theo}[\cite{KP}]
The moduli space $\MK{1,d}{X}$ is connected for any rational homogeneous space $X$ and any degree $d$.
\end{theo}

Note, in contrast, that for any $X$ in the tables of the introduction, we have the proposition.

\begin{prop}
\label{non-irred}
For $d\geq d(X)$, the moduli space $\MK{1,d}{X}$ is never irreducible.
\end{prop}

\begin{proo}
Indeed, consider the locally closed subset of $\MK{1,d}{X}$ where the degree on the unique elliptic irreducible component of the source of the map is 2. The dimension of this locally closed subset is $dc_1(X)+\dim X$ which is always strictly bigger than the expected dimension $dc_1(X)$ which is the dimension of the component containing irreducible curves for $d\geq d(X)$.
\end{proo}

However, there is a natural decomposition of this space into locally closed subsets $\MKE{1,d_\tau}{\tau}{X}$ parametrised by their combinatorial graph $\tau$ and combinatorial degree $d_\tau$ (see \cite{KM} or also \cite{KP}). Note that for stable maps of genus one, at most one irreducible component of the source map is an elliptic curve. We call this component (if it exists) the elliptic component of $\tau$. As a consequence of our results we obtain the following proposition.

\begin{prop}
If the degree $d_\tau$ is bigger that $d(X)$ on the elliptic component of $\tau$, then $\MKE{1,d_\tau}{\tau}{X}$ is irreducible.
\end{prop}

\begin{proo}
If there is no elliptic component, then the result follows from the fact (see \cite{BCMP}) that the moduli space of stable maps of genus 0 passing through at most 2 fixed point is irreducible for any homogeneous space $X$.

If there is an elliptic component, then claim that the natural forgetful map to the moduli space of elliptic curves is flat with irreducible fibers. Indeed, the fiber over $C$ is isomorphic to the product of ${\bf Hom}_d(C,X)$ with some moduli spaces of rational stable maps passing through at most 2 fixed points. These last moduli spaces are irreducible of the expected dimension by \cite[Corollary 3.3]{BCMP}. Now by dimension count, any irreducible component of $\MKE{d_\tau}{\tau}{X}$ has to dominate the moduli space of elliptic curves (the fibers have constant dimension 1 less than the expected dimension) thus the map is flat. The result follows by \cite[Corollaire 2.3.5.(\i\i)]{EGA}.
\end{proo}

%%%%%%%%%%%%%%%%%%%%%%%%%%%%%% Appendice %%%%%%%%%%%%%%%%%%%%%

\appendix
\section{The quadric of dimension 3}\label{sec:q3}

In this section, we prove Proposition \ref{prop:q3} for $d\geq 4$ using Bott-Samelson resolutions (the reader could see \cite{demazure} for more details on Bott-Samelson resolutions). This method was already used by the second author in \cite{arxiv} to study elliptic curves in spinorial varieties.

First we introduce some notation. Let $V$ be a 5-dimensional vector space endowed with a non degenerate quadratic form $Q$.
Then $X$ is the subvariety of $\p(V)$ consisting of isotropic lines in $V$.

For any isotropic flag $W_{\bullet}=(W_1,W_2)$ of $V$, we define a Bott-Sameslon resolution $\pi:\Xt\to X$ by $$\Xt=\Xt_{W_\bullet}=\{(V_1',V_2,V_1)\in \Q_3\times\G_Q(2,5)\times \Q_3\ \mid\ V'_1\subset W_2,\ V'_1\subset V_2,\ V_1\subset V_2\},$$ and $\pi(V_1',V_2,V_1)=V_1$. The map $\pi$ is birational, in particular it is an isomorphism over the open $B$-orbit of $X$, where $B$ is the Borel subgroup of $\operatorname{SO}(5)$ associated to $W_{\bullet}$. Moreover if we define $X_0=\{{\rm pt}\}$, $X_1=\{V'_1\in \Q_3\ \mid\ V'_1\subset W_2\}$ $X_2=\{(V_1',V_2)\in \Q_3\times\G_Q(2,5)\ \mid\ V_2\supset V'_1\subset W_2\}$ and $X_3=\Xt$,
there is sequence of $\p^1$-bundles
$$\Xt=X_3\to X_2\to X_1\to X_0=\{{\rm pt}\}.$$

The proof of Proposition \ref{prop:q3} uses the 3 following facts.
\begin{enumerate}[{Fact} 1]
\item Under some conditions $(*)_d$ on the class $\tilde{\alpha}$ of 1-cycles in $\Xt$ (see Definition \ref{defi:conditions}), the scheme  ${\bf Hom}_{\tilde{\alpha}}(C,\Xt)$ of morphisms from $C$ to $\Xt$ of class $\tilde{\alpha}$ is irreducible of dimension at most $3d-1$ with equality for a unique class $\tilde{\alpha}$).
\item For all $f\in{\bf Hom}_{d}(C,X)$ such that $f(C)$ is not contained in an isotropic projective line of $\p(V)$, we can choose a flag $W_\bullet$ such that $f$ lifts into a unique $\tilde{f}:C\to\Xt_{W_\bullet}$, such that $\tilde{\alpha}:=[\tilde{f}_*(C)]$ satisfies the conditions $(*)_d$. And, for each such $f$, the set of flags that we can choose is a subvariety (not closed) of $\operatorname{SO}(5)/B$ of dimension $3$.
\item  For $d\geq 4$ we prove that the scheme of morphisms from $C$ to $X$ of degree $d$ such that $f(C)$ is contained in an isotropic projective line of $\p(V)$ cannot be an irreducible component of ${\bf Hom}_{d}(C,X)$.
\end{enumerate}

Assuming these facts for the moment, let us prove Proposition \ref{prop:q3}.
Let $I$ be the set of couples $(W_\bullet,\tilde{f})\in\operatorname{SO}(5)/B\times{\bf Hom}_{(*)_d}(C,\Xt_{W_\bullet})$ where ${\bf Hom}_{(*)_d}(C,\Xt_{W_\bullet})$ is the scheme of morphisms $\tilde{f}$ from $C$ to $\Xt_{W_\bullet}$ such that $[\tilde{f}_*(C)]$ satisfies $(*)_d$. Call $p$ and $q$ the natural projections from $I$ to ${\bf Hom}_{d}(C,X)$ and to $\operatorname{SO}(5)/B$ respectively. Then, by the second and third facts, $p$ is dominant as soon as $d\geq 4$ and general fibers of $p$ have dimension~3.
Now, for any flag $W_\bullet$, the fiber $q^{-1}(W_\bullet)$ is the disjoint union, over the $\tilde{\alpha}$ satisfying conditions $(*)_d$, of the schemes ${\bf Hom}_{\tilde{\alpha}}(C,\Xt_{W_\bullet})$ of morphisms $\tilde{f}$ from $C$ to $\Xt_{W_\bullet}$ such that  $[\tilde{f}_*(C)]=\tilde{\alpha}$. By the first fact, exactly one of these latter schemes has maximal dimension $3d-1$ and moreover this scheme is irreducible, so that the fiber $q^{-1}(W_\bullet)$ is irreducible. Thus $I$ is a union of irreducible connected components of dimension at most $3d-1+=3d+3$ with equality for a unique component. The image in ${\bf Hom}_{d}(C,X)$ of these components are irreducible of dimension at most $3d+3-3=3d$ which is the expected dimension. Therefore all these images are contained in the closure of the image of the maximal one proving the irreducibility. Moreover the dimension of ${\bf Hom}_{d}(C,X)$ therefore $4+(3d-1)-3=3d$.\\

Let us prove the three facts one by one.
We begin with a general lemma.
\begin{lemm}
Let $\phi:X\to Y$ be a $\p^1$-fibration with a section $\sigma$. Let $\tilde{\alpha}$ be a class of 1-cycles in $X$. Denote by $T_\phi$ the relative tangent bundle and by $\xi$ the divisor $\sigma(Y)$. Suppose that $\tilde{\alpha}\cdot\xi\geq 0$ and that $\tilde{\alpha}\cdot(T_\phi-\xi)>0$.

Then, if ${\bf Hom}_{\phi_*\tilde{\alpha}}(C,Y)$ is irreducible, ${\bf Hom}_{\tilde{\alpha}}(C,X)$ is also irreducible, and $$\dim({\bf Hom}_{\tilde{\alpha}}(C,X))=\dim({\bf Hom}_{\phi_*\tilde{\alpha}}(C,Y))+\tilde{\alpha}\cdot T_\phi.$$
\end{lemm}

\begin{proo}
Let $E$ a rank two vector bundle on $Y$ such that $X=\p_Y(E)$. The section $\s$ is given by a surjection $E\to L$ where $L$ is an invertible sheaf on $Y$. Let $N$ be the kernel of this map.

We first study the fiber of the map $\MorC{\at}{X}\to\MorC{\vp_*\at}{Y}$ given by composition of morphisms over a morphism $f:C\to Y$. An element in this fiber  is given by a lift of $f$ \emph{i.e.} a surjective map $f^*E\to M$ of vector bundles on $C$ where $M$ is an invertible sheaf with $2\deg(M)-\deg(f^*E)=\at\cdot T_\vp$. An element in the fiber is therefore given by an invertible sheaf $M$ on $C$ with degree $d=\frac{1}{2}(\deg(f^*E)+\at\cdot T_\vp)$ and a surjective map in $\p(\Hom(f^*E,M))$.

Note that we have the equalities $\at\cdot\xi=\deg(M)-\deg(f^*N)$ and $\at\cdot(T_\vp-\xi)=\deg(M)-\deg(f^*L)$. We discuss two cases.

If $\at\cdot\xi>0$ then $\Hom(f^*E,M)$ is isomorphic to $\Hom(f^*N,M)\oplus\Hom(f^*L,M)$ and its dimension $\at\cdot T_\vp$ does not depend on $f$. For any invertible sheaf $M$, there exist surjections $f^*E\to M$ and they form an open subset of $\p\Hom(f^*E,M)\times\pic_{d}(C)$. Doing this construction in family as in \cite[Proposition 4]{fourier} we get a smooth fibration or relative dimension $\at\cdot T_\vp$ over $\MorC{\vp_*\at}{Y}$ and the result follows.

If $\at\cdot\xi=0$, then if $M\not\simeq f^*N$ we have $\Hom(f^*N,M)=0$ and any map $f^*E\to M$ factorises through $f^*L$ and is never surjective because $\deg(M)-\deg(f^*L)=\at\cdot(T_\vp-\xi)>0$. Therefore any pair $(M,p:f^*E\to M)$ of the fiber satisfies $M\simeq f^*N$. In that case, because of the equality $\deg(M)-\deg(f^*L)=\at\cdot(T_\vp-\xi)>0$, the sheaf $f^*E$ is isomorphic to $M\oplus f^*L$ and we have an isomorphism $\Hom(f^*E,M)\simeq\Hom(f^*N,M)\oplus\Hom(f^*L,M)$. The dimension of $\Hom(f^*E,M)$ is $\at\cdot T_\vp+1$ and does not depend on $f$. The fiber is therefore a non empty open subset of $\p\Hom(f^*E,M)$. We therefore again get a smooth fibration of relative dimension $\at\cdot T_\vp$ over $\MorC{\vp_*\at}{Y}$.
\end{proo}

Applying this lemma to the $\p^1$-fibrations $\phi_i:X_i\to X_{i-1}$ for $i\in\{1,2,3\}$ given above by a Bott-Samelson resolution $\Xt$ of $\Q_3$, we obtain the following corollary. We denote by $T_i$ the relative tangent space of $\phi_i$, we define three divisors $\xi_1$, $\xi_2$, $\xi_3$ of $\Xt$, which are given by natural sections of the $\phi_i$,  by
$$\begin{array}{c}
\xi_1=\{(V_1',V_2,V_1)\in\Xt\ \mid \ V_1'=W_1\}\\
 \xi_2=\{(V_1',V_2,V_1)\in\Xt\ \mid \ V_2=W_2\}\\
 \xi_3=\{(V_1',V_2,V_1)\in\Xt\ \mid \ V_1=V_1'\}.
 \end{array}
 $$

\begin{coro}
Let $\tilde{\alpha}$ be a class of 1-cycles in $\Xt$ satisfying for all $i\in\{1,2,3\}$, $\tilde{\alpha}\cdot\xi_i\geq 0$ and $\tilde{\alpha}\cdot(T_i-\xi_i)>0$.

Then ${\bf Hom}_{\tilde{\alpha}}(C,\Xt)$ is irreducible of dimension $\tilde{\alpha}\cdot(T_1+T_2+T_3)$.
\end{coro}

To obtain the first fact, we use the following result that can be deduced from \cite[Corollary 3.8 and Fact 3.7]{J.alg} by a short computation.
\begin{prop}

(\i) The relative tangent bundles $T_i$ are expressed in terms of the $\xi_i$ as follows:
$$T_1=2\xi_1,\ \ T_2=2\xi_2+\xi_1\ \textrm{and} \ T_3=2\xi_3+2\xi_2.$$
 In particular, we have $\tilde{\alpha}\cdot\xi_i\geq 0$ and $\tilde{\alpha}\cdot(T_i-\xi_i)>0$ for all $i\in\{1,2,3\}$ as soon as $\tilde{\alpha}\cdot\xi_1>0$, $\tilde{\alpha}\cdot\xi_2>0$ and $\tilde{\alpha}\cdot\xi_3\geq 0$, or, $\tilde{\alpha}\cdot\xi_1>0$, $\tilde{\alpha}\cdot\xi_2\geq 0$ and $\tilde{\alpha}\cdot\xi_3>0$

(\i\i) The pull-back $\pi^*\mathcal{O}_{\Q_3}(1)$ of the ample generator of the Picard group of $\Q_3$ equals $\xi_1+2\xi_2+\xi_3$. In particular, if $\tilde{\alpha}=\pi^*(\alpha)$ where $\alpha$ is the class of 1-cycles of degree $d$ in $\Q_3$, we have $d=\tilde{\alpha}\cdot(\xi_1+2\xi_2+\xi_3)$.
\end{prop}

\begin{defi}\label{defi:conditions}
We say that $\tilde{\alpha}$ satisfies conditions $(*)_d$ if
\begin{itemize}
\item $\tilde{\alpha}\cdot\xi_1>0$, $\tilde{\alpha}\cdot\xi_2>0$ and $\tilde{\alpha}\cdot\xi_3\geq 0$, or $\tilde{\alpha}\cdot\xi_1>0$, $\tilde{\alpha}\cdot\xi_2\geq 0$ and $\tilde{\alpha}\cdot\xi_3>0$;
\item $d=\tilde{\alpha}\cdot(\xi_1+2\xi_2+\xi_3)$.
\end{itemize}
\end{defi}

Remark that, under conditions $(*)_d$, we have $\tilde{\alpha}\cdot(T_1+T_2+T_3)=3d-\tilde{\alpha}\cdot(2\xi_2+\xi_3)$ and then the dimension of the scheme ${\bf Hom}_{\tilde{\alpha}}(C,\Xt)$ is at most $3d-1$ with equality if an only if $\tilde{\alpha}\cdot\xi_1=d-1$, $\tilde{\alpha}\cdot\xi_2=0$ and $\tilde{\alpha}\cdot\xi_3=1$. We deduce the first fact from above results.\\

To prove the second fact, fix a curve $f\in{\bf Hom}_{d}(C,\Q_3)$ such that $f(C)$ is not contained in an isotropic projective line of $\p(V)$. We begin by proving the following lemma. If $E$ is a vector subspace of $V$, we denote by $E^\perp$ the orthogonal subspace to $E$ in $V$ (with respect to $Q$). For any point $x$ in $\p(V)$, we denote by $\overline{x}$ the corresponding line in $V$.

\begin{lemm}\label{lemm:flagchoice}
There exist 3 points $x_0$, $x_1$, $x_2$ of $f(C)$ and an isotropic flag $W_\bullet=(W_1,W_2)$ such that $\overline{x_0}$ is not contained in $W_1^\perp$, $\overline{x_1}$ is contained in $W_1^\perp$ but not in $W_2$, and $\overline{x_2}$ is contained in $W_2$.
\end{lemm}

\begin{proo}
First remark that, because $f(C)\subset\Q_3$, the condition saying that $f(C)$ is not contained in an isotropic projective line of $\p(V)$ means that $f(C)$ is not contained in a projective line of $\p(V)$.

Then there exist 3 non-collinear points $x_0$, $x_1$, $x_2$ in $f(C)$. Let $E:=\overline{x_0}\oplus\overline{x_1}\oplus\overline{x_2}$ and $F:=\overline{x_1}\oplus\overline{x_2}$. We can assume that $Q\vert_F$ is non degenerate, in particular $V=F\oplus F^\perp$ and $Q$ is non-degenerate on $F^\perp$. Then there exists an isotropic line $W_1$ in $F^\perp\backslash E^\perp$. And we can define $W_2=W_1\oplus\overline{x_2}$, because $\overline{x_2}\not\subset F^\perp$.

It is obvious that $\overline{x_2}$ is contained in $W_2$. And we can also easily check that $\overline{x_0}$ is not contained in $W_1^\perp$ and that $\overline{x_1}$ is contained in $W_1^\perp$ but not in $W_2$.
\end{proo}

\begin{coro}
There exists a flag $W_\bullet$ such that $f$ lifts into a unique $\tilde{f}:C\to\Xt_{W_\bullet}$, and $[\tilde{f}_*(C)]$ satisfies the conditions $(*)_d$.
\end{coro}

\begin{proo}
Let $W_\bullet$ be an isotropic flag and $(x_0,x_1,x_2)\in f(C)^3$ as in Lemma \ref{lemm:flagchoice}. Denote by $B$ the Borel subgroup of $\operatorname{SO}(5)$ corresponding to $W_\bullet$. Then the open $B$-orbit $\Omega$ in $\Q_3$ is the set of points $x$ such that $\overline{x}$ is not contained in $W_1^\perp$. Recall that, over $\Omega$, the morphism $\pi:\Xt_{W_\bullet}\to\Q_3$ is an isomorphism. Thus, since $C$ is smooth and $f(C)$ intersects $\Omega$ (at least in $x_0$), $f$ lifts into a unique $\tilde{f}:C\to\Xt_{W_\bullet}$.

Moreover, we can compute that $\pi(\xi_1)=\{x\in\Q_3\ \mid \ \overline{x}\subset W_1^\perp\}$ and $\pi(\xi_2)=\pi(\xi_3)=\{x\in\Q_3\ \mid \ \overline{x}\subset W_2\}$. So, for any $i\in\{1,2,3\}$, $\pi^{-1}(\Omega)$ does not intersect $\xi_i$ so that $[\tilde{f}_*(C)]\cdot\xi_i\geq 0$. Also, the existence of $x_1$ implies that $[\tilde{f}_*(C)]\cdot\xi_1>0$ and the existence of $x_2$ implies that $[\tilde{f}_*(C)]\cdot\xi_2>0$ or $[\tilde{f}_*(C)]\cdot\xi_3>0$.
\end{proo}

Let us now explain why the set of flags that we can choose is of dimension $3$. In the proof of Lemma \ref{lemm:flagchoice}, we can note that, when $x_0$, $x_1$, $x_2$ are fixed in $f(C)^3$, we choose $W_1$ in a quadric of dimension~1 and $W_2$ is uniquely determined. Moreover, if $W_\bullet$ is fixed, then we have finitely many choice for $x_1$ and $x_2$ whereas the set of possible $x_0$ is one-dimensional. At the end, the dimension we are looking for is $3+1-1=3$.\\

Finally, let us prove the third fact.

Let $J$ be the scheme of morphisms from $C$ to $X$ of degree $d$ such that $f(C)$ is contained in an isotropic projective line of $\p(V)$. We have a natural projection from $J$ to the variety $\G_Q(2,5)$ of isotropic projective lines in $\p(V)$, which is of dimension~3. The fibers are isomorphic to the scheme ${\bf Hom}_{d}(C,\p^1)$ of morphisms of degree $d$ from $C$ to $\p^1$ and are of dimension $2d$. The dimension of $J$ is therefore $2d+3$. If $d\geq 4$, we have $2d+3<3d$ and then $J$ cannot be an irreducible component of ${\bf Hom}_{d}(C,X)$.

%%%%%%%%%%%%%%%%%%%%%%% Biblio %%%%%%%%%%%%%%%%%%%%%%

\bigskip\noindent
Boris {\sc Pasquier}, \\
{\it Institut de Math{\'e}matiques et de Mod{\'e}lisation de Montpellier}
Universit{\'e} Montpellier 2, CC 51, place Eug{\`e}ne Bataillon, 34095 Montpellier Cedex 5, France.

\noindent {\it email}: \texttt{boris.pasquier@math.univ-montp2.fr}.

\medskip\noindent
Nicolas {\sc Perrin}, \\
{\it Hausdorff Center for Mathematics,}
Universit{\"a}t Bonn, Villa Maria, Endenicher
Allee 62,
53115 Bonn, Germany and \\
{\it Institut de Math{\'e}matiques de Jussieu,}
Universit{\'e} Pierre et Marie Curie, Case 247, 4 place
Jussieu, 75252 Paris Cedex 05, France.

\noindent {\it email}: \texttt{nicolas.perrin@hcm.uni-bonn.de}.


\begin{thebibliography}{BCMP10}

\bibitem[Bal89]{ballico} Ballico, E., {\it On the Hilbert scheme of curves in a smooth quadric}.  Deformations of mathematical structures (\L{\'o}d\'z/Lublin, 1985/87),  127{\^a}132, Kluwer Acad. Publ., Dordrecht, 1989.

\bibitem[Bou54]{bou} Bourbaki, N., {\it Groupes et alg{\`e}bres de Lie}.
  Hermann 1954.

\bibitem[Bru87]{bruguieres} Brugui{\`e}res, A., {\it The scheme of morphisms from an elliptic curve to a Grassmannian.} Compositio Math. {\bf 63} (1987), no. 1, 15-40.

\bibitem[Buc03]{buch} Buch, A.S., {\it Quantum cohomology of Grassmannians}.  Compositio Math. {\bf 137} (2003), no. 2, 227-235.

\bibitem[BCMP10]{BCMP} Buch, A.S., Chaput, P.-E., Mihalcea, L.C., Perrin, N., {\it Finiteness of cominuscule quantum K-theory}, arXiv:1011.6658.

\bibitem[CP09a]{rational} Chaput, P.-E., Perrin, N., {\it Rationality of some Gromov-Witten varieties and application to quantum K-theory}, arXiv:0905.4394. To appear in Comm. in Contemp. Math.

\bibitem[CP09b]{adjoint} Chaput, P.-E., Perrin, N., {\it On the quantum cohomology of adjoint homogeneous spaces}, arXiv:0904.4824. To appear in Proc. of the London Math. Soc.

\bibitem[Dem74]{demazure} Demazure, M., {\it D{\'e}singularisation
des vari{\'e}t{\'e}s de Schubert g{\'e}n{\'e}ralis{\'e}es.} Collection
of articles dedicated to Henri Cartan on the occasion of his 70th birthday, I.
Ann. Sci. {\'E}cole Norm. Sup. (4) {\bf 7} (1974) 53-88.

\bibitem[Gro61]{grothendieck} Grothendieck, A., {\it Techniques de construction et th{\'e}or{\`e}mes d'existence en g{\'e}om{\'e}trie alg{\'e}brique. IV. Les sch{\'e}mas de Hilbert.} S{\'e}minaire Bourbaki, Vol. 6,  Exp. No. 221, 249-276, Soc. Math. France, Paris, 1995.

\bibitem[Gro65]{EGA} Grothendieck, A., {\it El{\'e}ments de g{\'e}om{\'e}trie alg{\'e}brique. IV. {\'E}tude locale des sch{\'e}mas et des morphismes de sch{\'e}mas. II.} Inst. Hautes {\~A}tudes Sci. Publ. Math. No. 24 1965.

\bibitem[GLP83]{GLP} Gruson, L., Lazarsfeld, R., Peskine, C., {\it On a theorem of Castelnuovo, and the equations defining space curves.}
Invent. Math. {\bf 72} (1983), no. 3, 491{\^a}506.

\bibitem[Har77]{hartshorne} Hartshorne, R., {\it Algebraic geometry}. Graduate Texts in Mathematics, No. 52. Springer-Verlag, New York-Heidelberg, 1977.

\bibitem[IM07]{marku} Iliev, A., Markushevich, D., {\it Parametrization of sing $\Theta$ for a Fano 3-fold of genus 7 by moduli of vector bundles}.  Asian J. Math. {\bf 11}  (2007),  no. 3, 427-458.

\bibitem[KP01]{KP} Kim, B., Pandharipande, R., {\it The connectedness of the moduli space of maps to homogeneous spaces.} Symplectic geometry and mirror symmetry (Seoul, 2000), 187-201, World Sci. Publ., River Edge, NJ, 2001.

\bibitem[Kle74]{kleiman} Kleiman, S.L., {\it The transversality of a general translate.} Compositio Math. {\bf 28} (1974), 287-297.

\bibitem[KM94]{KM} Kontsevich, M., Manin, Y., {\it Gromov-Witten classes, quantum cohomology, and enumerative geometry}.  Comm. Math. Phys. {\bf 164}  (1994),  no. 3, 525-562.

\bibitem[Mor79]{MO} Mori, S., {\it Projective manifolds with ample tangent bundles.}  Ann. of Math. (2) {\bf 110} (1979), no. 3, 593-606.

\bibitem[Ott88]{ottaviani} Ottaviani, G., {\it Spinor bundles on quadrics.}
Trans. Amer. Math. Soc. {\bf 307} (1988), no. 1, 301-316.

\bibitem[Per02]{fourier} Perrin, N., {\it Courbes rationnelles sur les vari{\'e}t{\'e}s homog{\`e}nes}. Ann. Inst. Fourier (Grenoble) {\bf 52} (2002),  no. 1, 105-132.

\bibitem[Per05]{J.alg} Perrin, N., {\it Rational curves on minuscule Schubert varieties}, J. Alg. {\bf 294} (2005), 431-462.

\bibitem[Per06]{arxiv} Perrin, N., {\it Courbes elliptiques sur la vari{\'e}t{\'e} spinorielle}, arXiv:math/0607260.

\bibitem[Per07]{compositio} Perrin, N., {\it Small resolutions of minuscule Schubert varieties},  Compos. Math. {\bf 143}  (2007),  no. 5, 1255-1312.

\end{thebibliography}
\end{document}